\documentclass[12pt]{article}
\usepackage{amsmath}
\usepackage{amssymb}
\usepackage{amsthm}
\usepackage{graphicx}
\usepackage{cancel}
\usepackage[normalem]{ulem}
\usepackage[usenames]{color}

\usepackage{hyperref}

\numberwithin{equation}{section}

\newtheorem{thm}{Theorem}[section]

\newtheorem{rem}[thm]{Remark}

\newtheorem{definition}[thm]{Definition}

\renewcommand{\qed}{\ifmmode$\Box$\else{\unskip\nobreak\hfil
	\penalty50\hskip1em\null\nobreak\hfil$\Box$
	\parfillskip=0pt\finalhyphendemerits=0\endgraf}\fi}

\newcommand{\eps}{\varepsilon}
\newcommand{\bfu}{{\boldsymbol u}}
\newcommand{\bfv}{{\boldsymbol v}}
\newcommand{\bfz}{{\boldsymbol z}}

\newcommand{\bfw}{{\boldsymbol w}}
\newcommand{\bfg}{{\boldsymbol g}}
\newcommand{\bfh}{{\boldsymbol h}}
\newcommand{\bfW}{{\boldsymbol W}}
\newcommand{\bfH}{{\boldsymbol H}}
\newcommand{\bfV}{{\boldsymbol V}}

\begin{document}

\vspace*{1cm}

\begin{center}
{\bf \large Mathematical Modeling of Biofilm Development}
\vspace{1cm}

{\sc Maria Gokieli}, {\sc Nobuyuki Kenmochi}

and

{\sc Marek Niezg\'odka}

Interdisciplinary Centre for Mathematical\
and Computational Modelling,\\
University of Warsaw,\\
Pawi\'nskiego 5a, 02-106 Warsaw, Poland
\end{center}
\vspace{1cm}

\noindent
{\bf Abstract.} We perform mathematical anaysis of the biofilm development process.
A~model describing biomass growth is proposed: 
 It arises from coupling three parabolic nonlinear equations:  a biomass equation with degenerate and 
singular diffusion, a nutrient tranport equation with a biomass-density
dependent diffusion, and an equation of   the Navier-Stokes type, 
describing the fluid flow in which the biofilm develops. 
This flow is subject to a biomass--density dependent obstacle.
The model is treated as a system of three inclusions, or variational inequalities;  
the third  one causes major difficulties for the system's
solvability. Our approach is based on the recent development of the theory on Navier-Stokes
variational inequalities. 
\renewcommand{\large}{\normalsize}
\renewcommand{\Large}{\normalsize}
\renewcommand{\LARGE}{\normalsize}

\section{Introduction}\label{sec:intro}
It is quite important for our furture to find clean and reproducible materials
and energy resources. In this connection, biomass has been noticed 
for the last thirty years.  Biomass growth is a process of aggregation of some living organisms
transported in fluids (liquids or gaz), usually sticking to the walls of the fluid container,
and thus influencing the flow itself.  It also involves nutrient 
transport and consumption. It can occur in air, water, soil penetrated by any fluid, blood. 
Only little is known about mathematical models of this mechanism. 
In particular, the process occurs in fluids, but
models coupled with hydrodynamics have been seldom analysed.

In~[8], such a biomass growth model coupled with fluid dynamics
has been proposed in the three dimensional space.  However, as far as we know, 
no theoretical results appeared in this context. 
The model assumed a sharp interface between the (solid) biomass and the liquid.
In the present paper we propose an analogous mathematical model of biomass growth dynamics
in a fluid, postulating, in place of a sharp interface, 
a thicker layer, considered as a
mixture of both phases --- just as in the weak formulation of a solid--liquid phase
transition.

For other formulations of biomass growth with taxis terms, 
see [7]. These formulations are not explicitly included 
in our formulation, but can be easily obtained by a modification.
\vspace{0.25cm}

 Let us recall in more detail the mathematical full model
proposed in~[8]. Let $\Omega \subset {\bf R}^3$ be a 
container in which biomass growth takes place. 
The process is described in terms of three unknown functions 
${\bfv}(x,t)$,  $w(x,t)$ and $u(x,t)$ which are respectively the velocity of
the fluid, the nutrient concentration and the biomass density at a point 
$x \in \Omega$ and time $t\geq 0$. They are governed by the following system:
\noindent
\begin{eqnarray*}
&(H_0) &{\bfv}_t+({\bfv}\cdot \nabla){\bfv} -\nu \Delta{\bfv}
   =-\frac 1\rho \nabla P,
~{\rm div~}{\bfv}=0,~{\rm in~}\{(x,t)|u(x,t)=0\},\\[0.2cm]
& &~~~~~~~{\rm where~}\rho ~{\rm is~the~constant~density~and~}P~{\rm is~the~ 
pressure~ in~ the~ fluid,}\\[0.2cm]
&(N_0)& w_t+{\bfv}\cdot \nabla w
     - {\rm div~}(d(u)\nabla w)=-f(w)u
~{\rm in~}\Omega,~t>0,\\
& &~~~~~~~{\rm where~}
    f(w)=\frac {k_1w}{k_2+w}~{\rm for~positive~constants~}k_1,~k_2,
    \\[0.2cm]
&(B_0)& u_t-\Delta d_1(u)+bu=f(w)u
~~{\rm in~}\Omega,~t>0,
\end{eqnarray*}
subject to suitable initial and boundary conditions. This model is derived
under the postulate that
the fluid cannot penetrate into the solid biomass ($u>0$),
the nutrient is convective by ${\bfv}\cdot \nabla w$ and diffusive with 
biomass-density dependent coefficient $d(u)$, and the diffusion of biomass is very slow near the interface
$u=0$, but very fast near the maximum density $u=u^*$.
 The function $f$ is the nutrient consumption term 
 and $b$ is a positive constant.

In this paper, we propose 
some relaxations and modifications into the above model, postulating that:
\begin{description}
\item{(i)} The biomass density $u(x,t)$ is non-negative and it has the finite maximum value $u^*$, i.e.\ $0\le u(x,t) \le u^*$. 
For some $\delta_0\in(0,u^*)$, which is fixed,
we postulate that the region of high density 
$\delta_0\le u(x,t) \le u^*$ 
is solid, and that of low density $0<u(x,t)<\delta_0$ 
is the interface layer between the solid biomass and the liquid. 
In such a layer, the behavior of $u$ may
correspond to the dynamics of 
planktonic biomass floating in the liquid, cf.~e.g.~[17]. 

This causes a biomass dependent constraint on the fluid's velocity.  
The constraint is written as:
  $$ |{\bfv}(x,t)| \le p_0(u^\varepsilon(x,t)), $$
where $p_0(r):(0,u^*]\to [0,\infty)$ is a $C^1$, non-negative and 
non-increasing function on $(0,u^*]$ such that (see Fig.1(i)):
 \begin{equation}\label{eq:1.4}
  \lim_{r \downarrow 0}p_0(r)=\infty,~~p'_0(r) <0, 
     ~\forall r \in(0,\delta_0),~~p_0(r)=0, ~\forall r \in 
    [\delta_0,u^*]; 
 \end{equation}
on the other hand  
$u^\varepsilon:=\rho_\varepsilon*u$ is the local spatial-average of $u(x,t)$ 
by means of the usual mollifier 
$\rho_\varepsilon(x)$ (see Section 2 for details).

\item{(ii)} The nutrient concentration $w(x,t)$ is non-negative 
and
has the threshold value $1$, i.e.\ $0\le w(x,t)\le 1$. Also, we suppose that 
there is no nutrient supply from the exterior. The diffusion coefficient
$d(u)$ depends on the biomass density $u$ and
  \begin{equation}\label{eq:1.2} c_d \le d(r) \le c'_d,~~|d(r_1)-d(r_2)| \le L(d)|r_1-r_2|,
   ~~\forall r_1,~r_2 \in {\bf R}, 
   \end{equation}
where $c_d,~c'_d$ and $L(d)$ are positive constants (see Fig.1(ii)). 
The function $f(w)u$, 
appearing in biomass density and nutrient transport equations, is called
the nutrient consumption, and in our model we suppose that
 \begin{equation}\label{eq:1.3}
  f(w) ~{\rm ~is ~of~}C^1~{\rm and~Lipschitz~ in~} w \in {\bf R},~
  f(0)\le 0~{\rm and~} f(1) \geq 0.
  \end{equation}
  
\item{(iii)} 
Biomass is diffusive 
(slowly near $u=0$, but fast near $u=u^*$), as well as convective 
by ${\bfv}\cdot \nabla u$. 
The degenerate diffusion term $d_1(u)$ is strictly increasing in $u\in [0,u^*)$ and
\begin{equation} \label{eq:1.1}
 d_1(0)=0,~~ \lim_{r \downarrow 0}\frac {d_1(r)}r =0,
     ~~\lim_{u \uparrow u^*}
      d_1(u)=\infty~~~~(\text{see~Fig.1(iii)}). 
      \end{equation}
      Note that we do not suppose $d_1$ to be continuous.

\end{description}

\begin{figure}
(i)
\includegraphics[width=4.5cm, height=3.8cm]{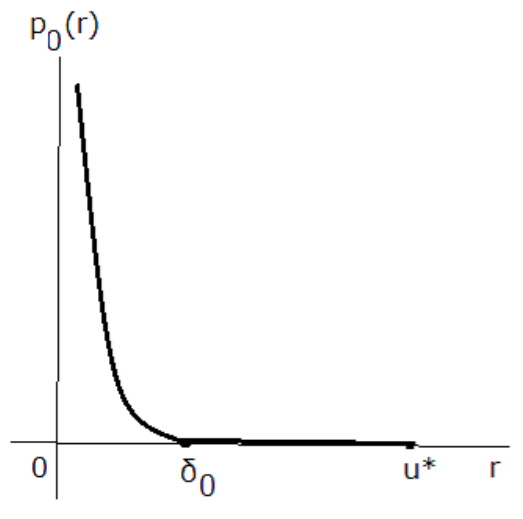}
(ii)
\includegraphics[width=4.5cm, height=3.8cm]{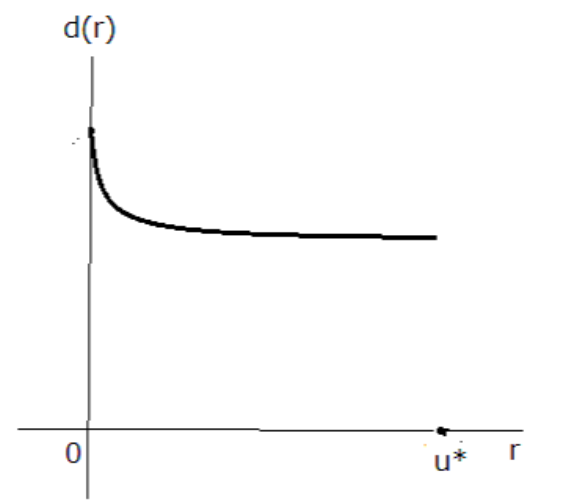}
(iii)
\includegraphics[width=4.5cm, height=3.8cm]{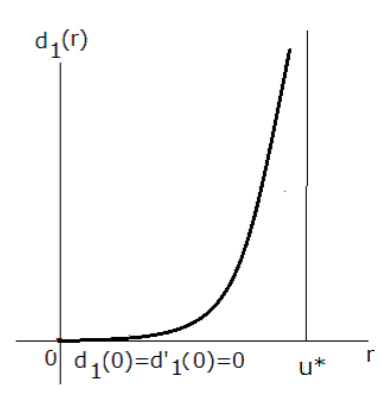}
\caption{Functions introduced in (i)--(iii): the obstacle function $p_0$, governing the flow velocity, the nutrient diffusivity~$d$, the biomass diffusivity $d_1$.}
\end{figure}

Now our relaxed/modified version for $\{(H_0),(N_0),(B_0)\}$ is described
as a system of three evolution equations --- one of them with a constraint ---
which is of the form:
 \begin{eqnarray*}
 &(H)^\varepsilon & 
 {\bfv}_t + ({\bfv}\cdot \nabla){\bfv}-\nu \Delta {\bfv} 
=	\bfg
   -\frac 1\rho \nabla P,~{\rm div}~{\bfv}=0,
   ~|{\bfv}|\le p_0(u^\varepsilon),~~~\text{in~}\Omega,~t>0,\\[0.2cm]
 &(N)& w_t +{\bfv}
     \cdot \nabla w -{\rm div}(d(u)\nabla w)  =  -f(w)u
     ~~{\rm in~}\Omega,~t>0,\\[0.2cm]
 &(B)& u_t +{\bfv}\cdot \nabla u
      -\Delta d_1(u) +b u = f(w)u~~{\rm in~}\Omega,~t>0.
 \end{eqnarray*}
The term $\bfg$ is an external force. As for boundary conditions, we take a standard Dirichlet boundary condition 
on the velocity $\bfv$,
--- 
which, without loss of generality, can be supposed homogenous
--- 
a homogenous Neumann boundary condition on the nutrient concentration $w$, 
and a mixed boundary condition on the biomass density $u$. 
The last means  a homogenous Neumann boundary condition for $u$ on all but some part
 of the boundary, $\Gamma_0\subset\partial\Omega$, which is supposed not to be touched
by the growing biomass: $u=0$ on $\Gamma_0$. 
We follow here [17], where
$\Gamma_0$ is the part of the  boundary through which the flow goes~in.
 
We have three main points in which this relaxed/modified 
model differs, formally, from $\{(H_0),(N_0),(B_0)\}$:
the convective term in $(B)$, 
the obstacle function $p_0$ in $(H)^\eps$,
and  the additional parameter $\eps$ 
---  actually two parameters, as  another one, $\delta_0$, is present in $p_0$.
All of these points are related to the planktonic layer introduced in (i).
The first one is its most natural consequence: 
the plankton is transported. 
The second one is related to the same assumption, 
and is also a mathematical tool
crucial for our treatment. Note that $(H_0)$ includes a constraint, 
meaning: no flow when $u>0$, free flow when $u=0$. This is a sharp interface model. 
The constraint in $(H)^\eps$, expressed in terms of 
$p_0$, is a blurred version of the previous one. 
The 'blurring' is governed by two parameters, $\eps$ and~$\delta_0$.
As a matter of fact,  we may reduce the number of parameters 
by taking $\delta_0= \delta_0(\eps)$ 
with $\delta_0(\eps)\downarrow 0$ as $\eps\downarrow 0$;
still,  as they are  independent, we leave both.
When $\varepsilon \downarrow 0$ and 
$\delta_0 \downarrow 0$ in  $\{(H)^\varepsilon,(N),\,
(B)\}$, we formally arrive at $\{(H_0),(N_0),(B_0)\}$. 
However, it seems quite difficult to carry out rigorously
this limit procedure. 

The main objective of this paper is to give an 
existence result for $\{(H)^\varepsilon,\,(N),\,(B)\}$, 
fixing parameters $\varepsilon>0$ and $\delta_0>0$.
The result is completely new and the model itself   
 reasonable from the biological point of view, despite the approximating
parameters. 

From the mathematical point of view,
$(H)^\varepsilon$ is going to be formulated in the solenoidal function space 
${\bfH}^1_{0,\sigma}(\Omega)$, $(N)$ and $(B)$ in the dual 
space of $H^1(\Omega)$.
Each problem  $(H)^\varepsilon$,
$(N)$ and $(B)$ is separately treated in the 
above-mentioned spaces (cf. [3, 5, 6, 9, 10]). However, 
the structure of~our system $\{(H)^\varepsilon,\,(N),\,(B)\}$ 
is extremely complicated because of 
its quasi-variational structure (cf. [11, 15]). 
The main difficulty for the analysis 
arises from this complexity of the couplings, especially the one  
in $(H)^\eps$, which appears via the nonlinear and unknown--dependent constraint.

The organization of this paper is as follows. In  section 2, we introduce the analytical framework.
In~sections~3, 4 and~5, we formulate each model apart:
the biomass density evolution, the nutrient transport and the flow 
governed by a Navier-Stokes variational inequality, respectively.
We~also give a  smooth approximation for each model  and prove its convergence.
Finally, in section~6, we formulate an approximate full system by  coupling 
these three models,  and prove existence of its solution
by the Schauder fixed-point argument. Then, we
construct a solution of our original problem
$\{(H)^\varepsilon,\,(N),\,(B)\}$ as a limit of approximate solutions,
making use of a recent important development on 
variational inequalities of the Navier-Stokes type, see [12].
Our main result is Theorem~6.2.

\section{Functional framework}\label{sec:framework}

\noindent
\emph{2.1  Functionals and their subdifferentials}\vspace{0.25cm}

\noindent
For a general (real) Banach space $X$ we denote by 
$X^*$ its dual. 
We denote by $|\cdot|_X$ and $|\cdot|_{X^*}$ the norms
in $X$ and $X^*$, and by $\langle \cdot, \cdot \rangle_{X^*,X}$
the  duality pairing between both spaces.

\medskip

Now, let $X$ be reflexive and consider a functional $\psi:\: X \to {\bf R}\cup \{\infty\}$. We say that:
\begin{description}
\item[$\psi$ is proper,] if $-\infty < \psi(z) \le \infty$ for all $z \in X$ and if it is not idetically $\infty$;
\item[$\psi$ is lower semi-continuous (l.s.c.) on $X$,] if 
$\liminf_{n\to \infty}\psi(z_n) \geq \psi(z)$ for any sequence $\{z_n\}$
converging to $z$ in $X$;
\item[$\psi$ is convex on $X$,] if $\psi(rz_1+(1-r)z_2) \le r\psi(z_1)
+(1-r)\psi(z_2)$ for all $z_1,\,z_2 \in X$ and $r \in [0,1].$\end{description}
For a proper, l.s.c.\ and convex function $\psi$ on $X$, the set 
 $$ D(\psi):=\{z \in X~|~\psi(z)< \infty\}$$
is called the \textbf{effective domain}. For each $z \in D(\psi)$
we consider a subset of $X^*$ 
 $$ \partial_{X^*,X}\psi(z):=\{z^* \in X^*~|~\langle z^*, v-z \rangle_{X^*,X} 
   \le \psi(v)-\psi(z), ~\forall v \in X\}, $$
which is called the \textbf{subdifferential} of $\psi$ at $z$; we put 
$\partial_{X^*,X}\psi(z)=\emptyset$ for $z \notin D(\psi)$.
If $X$ is a Hilbert space 
and it is identified with its dual, 
the subdifferential of a proper, l.s.c.~and convex function $\psi$ 
on $X$ is defined by using the inner product
$(\cdot,\cdot)_X$ in place of the duality $\langle \cdot,\cdot\rangle_{X^*,X}$
and the subdifferential at $z \in X$ is denoted by 
$\partial_X \psi(z)$:
 $$ \partial_{X}\psi(z):=\{ y \in X~|~(y, v-z )_{X} 
   \le \psi(v)-\psi(z), ~\forall v \in X\}. $$
For fundamental concepts and basic properties 
of subdifferentials we refer to [1, 4, 14]. 
\vspace{0.5cm}

\noindent
\emph{2.2 The domain}\vspace{0.25cm}

\noindent
Throughout this paper, we fix:
\begin{description} 
\item[$\Omega$,] a bounded domain in ${\bf R}^3$ with smooth boundary 
$\Gamma:=\partial \Omega$; 
\item[$\Gamma_0$,] 
a compact subset of $\Gamma$, having positive surface measure;
\item[$T$,] which is an arbitrary positive real number, and we denote $Q=\Omega\times[0,T]$.
\end{description}
\vspace{0.25cm}

\noindent
\emph{2.3 Function Spaces}\vspace{0.25cm}

\noindent
We set up:
$$ H:=L^2(\Omega),~~V:=H^1(\Omega).$$
The norms $|\cdot|_H$ and $|\cdot|_V$ are defined as usual. 
Next, 
denote by $V_0$ the space
 $$ V_0:=\{z \in V~|~ z=0~\text{ a.e.~on~}\Gamma_0 
\}, \quad
 \text{ with the norm }
   |z|_{V_0}:=|\nabla z|_H.$$
The condition $z=0$ above is understood in the sense of trace. We assume always that the dual spaces $V^*$ and $V^*_0$ 
are equipped with the dual norms of $V$ and $V_0$, respectively. 
By identifying $H$ with its dual space, we have
\begin{equation}\label{eq:2.1}
V\subset H \subset V^*,~~V_0 \subset H\subset V_0^*,
    ~~{\rm with~compact~embeddings}; 
\end{equation}
throughout this paper, we fix a positive
constant $c_0$ such that
\begin{equation}\label{eq:2.2}
|z|_H\le c_0 |z|_V~\forall z \in V,~~~
|z|_H\le c_0 |z|_{V_0}~\forall z \in V_0,~~~
|z|_{V^*}\le c_0 |z|_H,~|z|_{V_0^*}\le c_0 |z|_H~\forall z \in H.
\end{equation}
For simplicity of notation, 
the inner product $(\cdot,\cdot)_H$ in $H$, the dualities
$ \langle \cdot,\cdot \rangle_{V^*,V}$ and $\langle \cdot,\cdot 
\rangle_{V_0^*,V_0} $ are denoted  by $(\cdot,\cdot)$,  
$\langle \cdot,\cdot \rangle$ and $\langle \cdot,\cdot \rangle_0$, respectively. 

\medskip

The duality mapping $F_0$ from $V_0$ onto $V_0^*$ is characterized by 
\begin{equation}\label{eq:2.3}
 \langle F_0z_1,z_2 \rangle_0 =\int_\Omega \nabla z_1(x)\cdot \nabla z_2(x)dx
=: (F_0 z_1, F_0 z_2)_* 
, ~~\forall z_1,~z_2 \in V_0, 
\end{equation}
where the first equality defines $F_0$ and the second the induced inner product in $V_0^*$, denoted by $(\cdot, \cdot)_*$.
From the definition of $F_0$ and $V_0$, it follows (cf. [14; $\S 1$]) that formally 
\vspace{-0.24cm}
\begin{equation}\label{eq:2.4}
F_0v=-\Delta v ~{\rm in~}\Omega~{\rm in~the~sense~of~distributions},~
  v=0~{\rm on~}\Gamma_0,~\frac {\partial v}{\partial n}=0~{\rm on ~}
\Gamma-\Gamma_0.
\end{equation}

Next, we consider  solenoidal function spaces. Let
 $$\bf{\mathcal D}_\sigma:=\{{\bfz} \in~C^\infty_0(\Omega)^3~|~{\rm div}\bfz
=0~{\rm in}~\Omega\},$$
 $${\bfV}_\sigma= {\rm the~ closure~of~}{\bf{\mathcal D}}_\sigma~{\rm in~}
 H_0^1(\Omega)^3,~~{\bfH}_\sigma~={\rm ~ the~ closure~ of~} 
{\bf{\mathcal D}}_\sigma~{\rm in~}
L^2(\Omega)^3.$$ 
In these spaces the norms are given by
  $$ |{\bfz}|_{{\bfH}_\sigma}:
   =\left(\sum_{k=1}^3 \int_\Omega |z^{(k)}(x)|^2dx \right)^{\frac 12},
  ~~\forall {\bfz}=(z^{(1)},z^{(2)}, z^{(3)}) \in {\bfH}_\sigma,$$
  $$\text{and}~~ |{\bfz}|_{{\bfV}_\sigma}:
   =\left(\sum_{k=1}^3 \int_\Omega |\nabla z^{(k)}(x)|^2dx \right)^{\frac 12},
  ~~\forall {\bfz}=(z^{(1)},z^{(2)}, z^{(3)}) \in {\bfV}_\sigma.$$
Note that ${\bfH}_\sigma$ is a Hilbert space and by identifying it with
its dual, we have
 \begin{equation}\label{eq:2.5}
 {\bfV}_\sigma \subset {\bfH}_\sigma \subset {\bfV}^*_\sigma~~
   {\rm with~compact~embeddings}. 
 \end{equation}
We write $(\cdot,\cdot)_\sigma$ for the inner product in ${\bfH}_\sigma$ and $\langle \cdot,\cdot\rangle_\sigma$ for duality between 
${\bfV}^*_\sigma$ and ${\bfV}_\sigma$. \vspace{0.25cm}

\noindent
{\bf Remark 2.1.} We mean by $H \subset V^*_0$,  $H \subset V^*$, and
${\bfH}_\sigma \subset {\bfV}^*_\sigma$ 
in \eqref{eq:2.1} and \eqref{eq:2.5} that
$ \langle u, z \rangle_0=(u,z)$ for all $u \in H,~z\in V_0$ and 
$\langle u, z \rangle=(u,z)$ for all $u \in H, ~z\in V$ as well as
$ \langle {\bfu}, {\bfz} \rangle_\sigma=({\bfu},{\bfz})_\sigma$
for all ${\bfu} \in {\bfH}_\sigma, ~{\bfz}\in {\bfV}_\sigma$. 
\vspace{0.35cm}

\noindent
{\bf Remark 2.2.} 
If  $\bfv \in \bfV_\sigma$, then $\bfv = {\bf 0}$ 
on~$\partial \Omega$ and
${\bfv}\cdot \nabla z={\rm div}(z{\bfv})$ for all $z \in V$. 

\vspace{0.25cm}

\noindent
\emph{2.4 Space averaging}\vspace{0.25cm}

\noindent
Given $\mu \in (0,1]$, a function $u \in H$ and any smooth function $\gamma$ on ${\bf R}^3$, 
we denote by $\rho_\mu * (\gamma u)$ the convolution of the usual
mollifier 
  $$\rho_\mu(x):=\left\{
           \begin{array}{l}\displaystyle
           \frac{1}{N_\mu}
            {\rm exp}\left(-\frac 1{{\mu}^2-|x|^2}\right)~~{\rm if~}|x|< \mu,\\[0.12cm]
             0,~~~~~~~~~~~~~~~~~{\rm otherwise},
            \end{array}\right. \quad
            N_\mu = \int_\Omega \exp\left(\frac 1{{\mu}^2-|x|^2}\right)dx,$$
and function $\gamma(x)u(x)$, namely
  $$ [\rho_\mu*(\gamma u)](x):=\int_{{\bf R}^3}\rho_\mu(x-y)\gamma(y)
     \tilde u(y) dy,~~\forall x \in \Omega,
           $$
where $\tilde u$ denotes the extension of $u$
onto ${\bf R}^3$ by $0$. 
Noting here that 
 $$ [\rho_\mu*(\gamma u)](x)
  =\int_\Omega \rho_\mu(x-y)\gamma(y)u(y)dy
  =(u, \gamma \rho_\mu(x-\cdot)),$$
we see that, in the case when $\gamma =0$ on $\Gamma_0$
\begin{equation}\label{eq:2.6} 
|\rho_\mu*(\gamma u)|_{C(\overline\Omega)}
    \le \left( \sup_{x \in \Omega}
     |\gamma \rho_\mu(x-\cdot)|_{V_0}\right) 
      |u|_{V^*_0} \le c_0\left( \sup_{x \in \Omega}
     |\gamma \rho_\mu(x-\cdot)|_{V_0}\right) |u|_H,
\end{equation}
and in the case when $\gamma \equiv 1$ 
\begin{equation}\label{eq:2.7} 
|\rho_\mu* u|_{C(\overline\Omega)}
    \le \left( \sup_{x \in \Omega}
     |\rho_\mu(x-\cdot)|_{V}\right) 
      |u|_{V^*} \le c_0\left( \sup_{x \in \Omega}
     |\rho_\mu(x-\cdot)|_{V}\right) |u|_H,
\end{equation} 
see \eqref{eq:2.2}. Similarly, 
if $u \in W^{1,2}(0,T;V^*_0)$ and $\gamma=0$ on $\Gamma_0$, then
\begin{equation}\label{eq:2.8}
|\rho_\mu*(\gamma u)|_{W^{1,2}(0,T;C(\overline\Omega))}
    \le \left( \sup_{x \in \Omega}
     |\gamma\rho_\mu(x-\cdot)|_{V_0}\right) 
     |u|_{W^{1,2}(0,T;V^*_0)},
\end{equation}
and if $u \in W^{1,2}(0,T;V^*)$, then 
\begin{equation}\label{eq:2.9}
|\rho_\mu* u|_{W^{1,2}(0,T;C(\overline\Omega))}
    \le \left( \sup_{x \in \Omega}
     |\rho_\mu(x-\cdot)|_{V}\right) 
     |u|_{W^{1,2}(0,T;V^*)}.
\end{equation}

\section{Biomass growth inclusion and its approximation}\label{sec:biomass}
In order to describe the degenerate and singular diffusion for biomass
density we use a non-negative, proper, l.s.c. and convex function 
$\hat \beta(\cdot)$ on ${\bf R}$ given by:
  $$ \hat \beta(r):=\left \{
       \begin{array}{ll}
  \displaystyle{ \int_0^r d_1(s)ds,}&~~~~~~{\rm for~}r \in [0,u^*],\\[0.5cm]
         \infty,&~~~~~~{\rm otherwise},
       \end{array} \right. $$
where $d_1$ is the function introduced in (i) in the
introduction, satisfying \eqref{eq:1.1}.
Its subdifferential $\beta:=\partial \hat \beta$ in ${\bf R}$ is
equal to $d_1$ except on a countable set, where $d_1$ is not necessarily continuous.
In these points of discontinuity, it is given by  $[d^{-}_1(r),d^{+}_1(r)]$, where $d^{-}_1(r):=\lim_{s\uparrow r}d_1(s)$ and $d^{+}_1(r):=
\lim_{s\downarrow r}d_1(s)$ for $r \in (0,u^*)$, if $r\in (0,u^*)$, 
also, $\beta(0) = (-\infty,0]$ and $\beta(r) = \emptyset$ for $r <0$ or $r \geq u^*$.
Clearly, $D(\beta)=[0,u^*)$, 
$d_1(r)\in\beta(r)$ for any $r\in [0, u^*)$,
$R(\beta)={\bf R}$ and $\beta$ is strictly monotone in ${\bf R}$ 
(see Fig.1(iii)). 

Now, we define the function $\varphi$ on $V^*_0$ by 
   $$\varphi(z):=\left \{
          \begin{array}{ll}
  \displaystyle {\int_\Omega \hat\beta(z(x))dx,}
        &~~{\rm if~}z \in H~{\rm and~}\hat\beta(z) \in L^1(\Omega),\\[0.5cm]
  \infty, &~~{\rm otherwise.}
   \end{array}\right. $$
Clearly, $\varphi(\cdot)$ is non-negative, proper and
convex on $V^*_0$ with $D(\varphi)$ included in the subset 
$\{z \in H~|~ 0\le z \le u^*~\text{a.e.\  on  }\Omega\}$. 
It follows that $\varphi$ is l.s.c.\ on~$V_0^*$. 
Hence any level
set of $\varphi(\cdot)$ is compact in~$V^*_0$.
We denote by $\partial_*\varphi(\cdot)$ the subdifferential of 
$\varphi(\cdot)$ in $V^*_0$, namely
$$ \partial_*\varphi(z):= \partial_{V_0^*}\varphi(z)=\{z^*\in V^*_0~|~(z^*,v-z)_* \le \varphi(v)-\varphi(z),
     \forall v \in V^*_0\}.~~\forall z \in D(\varphi).$$
Then we know (cf. [5, 6]) that
\begin{equation}\label{eq:3.2} 
\partial_* \varphi(v)=\{F_0\tilde v~|~\tilde v \in V,\tilde v 
      \in \beta(v)~{\rm a.e.~on~}\Omega \},~~\forall v \in 
      D(\partial_*\varphi)~(\subset H). 
\end{equation}
Let $g\in L^2(0,T; V^*_0)$ and $u_0 \in D(\varphi)$.
We denote by 
$CP(\varphi; g,u_0)$ the Cauchy problem
 $$ u'(t)+\partial_*\varphi(u(t)) \ni g(t) ~~{\rm in~}V^*_0~{\rm for~
a.e.~}t \in [0,T],~~u(0)=u_0. $$
By the general theory of evolution equations (cf. Appendix I) this Cauchy problem
admits one and only one solution $u$ such that $u \in W^{1,2}(0,T;V^*_0)$ and
$t \to \varphi(u(t))$ is absolutely continuous on $[0,T]$.  
The following convergence result will be used later on.\vspace{0.5cm}

\noindent
{\bf Lemma 3.1.} {\it Let $u_0 \in H$ with $u_0 \in D(\varphi)$ 
and $\{g_n\}$ be a sequence in $L^2(0,T;V^*_0)$ such that $g_n \to g$ weakly
in $L^2(0,T;V^*_0)$ as $n \to \infty$. Then, the solution $u_n$ of 
$CP(\varphi;g_n,u_0)$ converges to the solution $u$ of $CP(\varphi; g,u_0)$
in $C([0,T];V^*_0)\cap L^2(Q)$ and weakly in $W^{1,2}(0,T;V^*_0)$.}
\vspace{0.25cm}

\noindent
{\bf Proof.} The convergences 
$u_n \to u$ weakly in $W^{1,2}(0,T;V^*_0)$ and strongly
in $C([0,T];V^*_0)$ are obtained by Proposition II of the Appendix 
(note that $D(\varphi)$ is compact in $V_0^*$). We show 
below the convergence in $L^2(Q)$. Taking the  
difference of two inclusions for $u_n$ and $u$, we have
by \eqref{eq:3.2}
 $$ u'_n-u'+F_0(\tilde u_n-\tilde u)=g_n- g
   ~~{\rm in ~}V^*_0,$$
where $\tilde u_n \in L^2(0,T;V_0)$ 
with $\tilde u_n \in \beta(u_n)$  a.e.~on $Q$ and $\tilde u \in L^2(0,T;V_0)$
with $\tilde u \in \beta(u)$ a.e.~on $Q$. Now, take the inner product between
both sides of the above relation and $u_n-u$ in $V^*_0$ to obtain
 $$ \frac 12 \frac d{dt}|u_n(t)-u(t)|^2_{V^*_0} + \int_\Omega (\tilde u_n(t)-
  \tilde u(t))(u_n(t)-u(t))dx \le (g_n(t)-g(t),u_n(t)-u(t))_*$$
for a.e.~$t \in [0,T]$.
Integrating this inequality in time over $[0,t]$ yields
 $$ \frac 12 |u_n(t)-u(t)|^2_{V^*_0}+\int_0^t\int_\Omega 
   (\tilde u_n-\tilde u)(u_n-u)dxd\tau \le \int_0^t(g_n-g,
    u_n-u)_*d\tau$$
for all $t \in [0,T]$,
whence, by monotonicity,
 $$ \lim_{n\to \infty}\int_Q (\tilde u_n-\tilde u)(u_n-u)dxdt
    =0. $$
We derive from this convergence that $u_n \to u$ in $L^2(\Omega)$. 
In fact, by the strict monotonicity of $\beta$ and $0 \in \beta(0)$, 
for any small $\delta>0$ there is
a constant $C_\delta \in (0,1)$ such that 
  $$  \tilde r_1-\tilde r_2 \geq C_\delta ~{\rm if~} r_1-r_2 \geq \delta,~
  r_1, r_2 \in D(\beta),~\tilde r_1 \in \beta(r_1)~{\rm and~}
     ~\tilde r_2 \in \beta(r_2). $$
Hence, putting $E_{n,\delta}:=\{(x,t) \in Q~|~|u_n(x,t)-u(x,t)|\geq \delta\}$,
we observe that
 $$ \begin{array}{l}
 \displaystyle{C_\delta \int_Q|u_n-u| dxdt = C_\delta \int_{E_{n,\delta}}
    |u_n-u| dxdt +C_\delta \int_{Q-E_{n,\delta}}|u_n-u| dxdt} \\[0.3cm]
 \displaystyle{\le  \int_Q (\tilde u_n-\tilde u)(u_n-u)dxdt 
     + \delta C_\delta T|\Omega|,}
    \end{array} $$
where $|\Omega|$ denotes the volume of $\Omega$.
Accordingly, $ \limsup_{n\to \infty}\int_Q |u_n-u|dxdt \le 
\delta T|\Omega|$. Since $\delta>0$ is arbitrary and $0\le u_n\le u^*$ a.e.
on $Q$, we have $u_n \to u$ in $L^2(Q)$.
\hfill $\Box$
\vspace{0.5cm}

\noindent
With the operator $\partial_*\varphi$, the biomass growth 
equation $(B)$ with formal boundary condition 
$u=0~{\rm on~}\Gamma_0\times(0,T)$ and 
$\frac{\partial u}{\partial n}=0~{\rm on~}(\Gamma-\Gamma_0)\times(0,T)$ 
(cf.~\eqref{eq:2.4}),
is reformulated as the Cauchy problem:
\begin{equation}\label{eq:3.3}
  (B;w,{\bfv};u_0)~~\left\{
          \begin{array}{l}
   \displaystyle{ u'(t)+ \partial_*\varphi(u(t))+{\bfv}(t)\cdot \nabla u(t) 
   +bu(t)
    \ni f(w(t))u(t)}
\quad\text{ in~}V_0^*,
\\[0.12cm]
    \displaystyle{u(0)=u_0,} 
      \end{array} \right.
\end{equation}
where $w$, $\bfv$, $u_0$ are given. More precisely, we have the following definition of solution.
\vspace{0.25cm}

\noindent
{\bf Definition 3.1.} {\it Let $w \in L^2(0,T;V)\cap L^\infty(Q)$, ${\bfv}\in 
L^2(0,T;{\bfV}_\sigma)$ and $u_0 \in H$ with $\hat\beta(u_0) \in L^1(\Omega)$.
Then, a function $u:[0,T]
\to V^*_0$ is called a solution to $(B;w,{\bfv},u_0)$, 
if $u \in W^{1,2}(0,T;V^*_0)$, $0\le u\le u^*$ a.e.~on~$Q$,
and for a.e.~$t \in (0,T)$, \eqref{eq:3.3} is satisfied.\\
    }

\noindent
Note that $\hat \beta(u_0) \in L^1(\Omega)$ implies $0\le u_0 \le u^*$ a.e.~on $\Omega$.
So as to be explicit for the sense of \eqref{eq:3.3}, we note that
on account of \eqref{eq:3.2} and Remarks 2.1,~2.2, 
the solution 
$u$ of   $(B;w,{\bfv},u_0)$ 
satisfies
the following variational equality: there is 
$\tilde u: [0,T] \to V_0$ 
such that
\begin{multline*}
~~~~~~~~~~~~~~~~~~~~~~~~~~~~~~~~~~~~~~~%
\tilde u(x,t) \in \beta(u(x,t))~~{\rm a.e.~on~}Q,\\
 \langle u'(t), z\rangle_0 + \int_\Omega \nabla \tilde u(t)\cdot \nabla z dx
    - \int_\Omega u(t){\bfv}(t)\cdot\nabla z dx +(b u(t),z)\\
=( f(w(t))u(t),z),~~\forall z \in V_0,~\text{ for~a.e.~}t 
    \in (0,T).
\end{multline*}

\noindent    
In order to solve $(B;w,{\bfv};u_0)$, we
approximate it by the following problem
including a real positive parameter $\mu \downarrow 0$:
 $$(B;w,{\bfv},u_0)_\mu ~~\left\{
          \begin{array}{l}
    \displaystyle{ u'(t)+ \partial_*\varphi(u(t))
      +{\bfv}(t)\cdot \nabla[\rho_\mu*(\gamma_\mu u(t))] +bu(t)}\\[0.2cm]
    ~~~~~~~~~~\ni f(\rho_\mu*w(t))u(t)
    ~~~{\rm in~}V^*_0~{\rm for~a.e.}~t\in [0,T],\\[0.2cm]
    \displaystyle{u(0)=u_0,} 
      \end{array}
   \right.$$
where $\{\gamma_\mu(\cdot)\}_{\mu \in (0,1]}$ 
is a family of smooth functions on  $\boldsymbol{R}^3$
such that
\begin{equation}\label{eq:3.4}\gamma_\mu(y) \ \ \left \{
 \begin{array}{l}
 =0,~~~~~~~{\rm if~}{\rm dist}(y,\Gamma_0)
 \le \frac 12 \mu,\\[0.2cm]
 \in [0,1],~~~{\rm if~}\frac 12 \mu
 < {\rm dist}(y,\Gamma_0) < \mu,\\[0.2cm]
 =1,~~~~~~~{\rm if~ dist}(y,\Gamma_0)\geq \mu,
 \end{array} \right.
 \end{equation}
for all $\mu \in (0,1]$ and $\gamma_\mu(\cdot)$ is continuous in $C(\bar{\Omega})$ with respect to $\mu\in(0,1]$.
We have $0\le \gamma_\mu(y)\le 1$ and $\gamma_\mu(y)\to 1$ for any $y\in\Omega$ as $\mu\downarrow 0$. 
\bigskip

\noindent
{\bf Remark 3.1. } When $\mu=0$, 
$(B;w,{\bfv}, u_0)_\mu = (B;w,{\bfv}, u_0)$.\vspace{0.5cm}

\noindent
{\bf Proposition 3.1.} {\it Assume that \eqref{eq:1.1} holds 
and let $\mu \in (0,1]$.
Let ${\bfv}$ and $w$ be given functions such that 
\begin{equation}\label{eq:3.5} 
\begin{array}{c}
 {\bfv} \in L^2(0,T;{\bfV}_\sigma)\cap L^\infty(0,T;{\bfH}_\sigma),\\[0.3cm]
  w \in W^{1,2}(0,T;V^*)\cap L^2(0,T;V),~
   0\le w\le 1 ~{\it a.e.~on~}Q. 
    \end{array}
    \end{equation}
Also, let $u_0 \in H$ be such that $\hat \beta(u_0) \in L^1(\Omega) $. 
Then, there exists one and only one
solution~$u$ to $(B;w,{\bfv},u_0)_\mu $. This solution is such that
$t \to |\hat\beta(u(t))|_{L^1(\Omega)}$ is absolutely continuous 
on $[0,T]$.
Moreover, there is a non-negative, bounded and non-decreasing function 
$B_0(\cdot)$ on 
$[0,\infty)\times [0,\infty)$, independent of the
parameter $\mu \in (0,1]$,
such that 
\begin{equation}\label{eq:3.6} 
|u|^2_{W^{1,2}(0,T;V^*_0)}+ 
      \sup_{t \in [0,T]}|\hat \beta(u(t))|_{L^1(\Omega)}
     \le B_0\left(|{\bfv}|_{L^2(0,T;{\bfH}_\sigma)}, 
|\hat \beta(u_0)|_{L^1(\Omega)}\right).  
\end{equation}
}

For the proof of Proposition 3.1 we prepare two
lemmas. \vspace{0.5cm}

\noindent
{\bf Lemma 3.2.} {\it Assuming  \eqref{eq:3.5} we have, for all $v \in V_0$ and $t \in [0,T]$,
  \begin{eqnarray*}
   |f(\rho_\mu*w(t))v|_{V_0}
  &\le& 3(c_0+1)\left\{L(f)|\Omega|^{\frac 12}
   \sup_{x \in \overline \Omega}
 |\rho_\mu(x-\cdot)|_V+\max_{0\le r \le 1}f(r)\right\} |v|_{V_0}\\
  &=:& M^\mu_1|v|_{V_0},
  \end{eqnarray*}  
 where $L(f)$ is the
Lipschitz constant of $f$ and $c_0$ is the 
constant from \eqref{eq:2.2}.  }\vspace{0.2cm}

\noindent
{\bf Proof.} First we note that
\begin{equation}\label{eq:3.7}
\begin{array}{ll}
 \displaystyle{
  |\nabla[\rho_\mu*w](x,t)|^2}
    &\displaystyle{= \sum_{i=1}^3|[\rho_{\mu,x_i}*w](x,t)|^2
 =\sum_{i=1}^3\left|\int_\Omega
    \rho_{\mu,x_i}(x-y)w(y,t)dy \right|^2} 
    \\[0.3cm]
 &\displaystyle{\le \sum_{i=1}^3\left(\int_\Omega
    |\rho_{\mu,x_i}(x-y)|dy \right)^2
    \le |\Omega|
    |\rho_\mu(x-\cdot)|_V^2.}
 \end{array} 
 \end{equation}
By \eqref{eq:3.7},
 \begin{eqnarray*}
 |f(\rho_\mu*w)v|^2_{V_0}
 &=& \int_\Omega |\nabla[f(\rho_\mu*w)v]|^2dx
\\
  &\le& L(f)^2|\nabla (\rho_\mu*w)|
 ^2_{C(\overline \Omega)}\int_\Omega|v|^2dx
 +\max_{0\le r \le 1}f(r)^2\int_\Omega
 |\nabla v|^2dx \\
 &\le& \left\{L(f)^2|\Omega|\sup_{x\in \overline\Omega}
 |\rho_\mu(x-\cdot)|_V^2
 +\max_{0\le r \le 1}f(r)^2\right\}(c_0^2+1)
 |v|_{V_0}^2.
\end{eqnarray*}
Thus the required inequality is obtained. \hfill $\Box$ 
\vspace{0.5cm}

\noindent
{\bf Lemma 3.3.} {\it Assuming \eqref{eq:3.5} we have, for all $z \in H$ and $ t \in [0,T]$:
  \begin{equation}\label{eq:3.8} 
  |{\bfv}(t)\cdot \nabla [\rho_\mu*
  (\gamma_\mu z)]|_{V^*_0}
  \le M_2^\mu 
 |z|_{V^*_0}|{\bfv}(t)|_{{\bfH}_\sigma},~~ 
 | f(\rho_\mu*w(t))z -bz|_{V^*_0} 
     \le M_3^\mu |z|_{V^*_0}, 
  \end{equation}

where 
$M^\mu_2 :=\sup_{x \in\overline\Omega}|
  \gamma_\mu \rho_\mu (x-\cdot)|_{V_0}$ and 
$M^\mu_3:= M^\mu_1+b$.  
}\vspace{0.5cm}

\noindent
{\bf Proof.} For any $z \in H$ we have by \eqref{eq:2.6} and Remarks~2.1, 2.2: 
    \begin{eqnarray*}
  & & |{\bfv}(t)\cdot \nabla [\rho_\mu*
  (\gamma_\mu z)]|_{V^*_0} \\
   &=&\sup_{ v\in V_0,|v|_{V_0}\le 1} 
   \langle {\bfv}(t)\cdot 
     \nabla [\rho_\mu*(\gamma_\mu z)], 
     v\rangle_0
   = \sup_{ v\in V_0,|v|_{V_0}\le 1}
  \int_\Omega {\rm div}[\rho_\mu*(\gamma_\mu
     z){\bfv}(t)]v dx \\
  &=& \sup_{ v\in V_0,|v|_{V_0}\le 1}
 \int_\Omega\left\{-[\rho_\mu*(\gamma_\mu z)]
 {\bfv}(t)\cdot \nabla v dx\right\}
  \le 
      |\rho_\mu*(\gamma_\mu z)|_{C(\overline \Omega)}
      |{\bfv}(t)|_{{\bfH}_\sigma}\\
 & \le& \left(\sup_{x \in \overline \Omega}
 |\gamma_\mu\rho_\mu  
  (x-\cdot)|_{V_0} \right)|z|_{V^*_0}|{\bfv}(t)|_{{\bfH}_\sigma} =
 M_2^\mu  |z|_{V^*_0}|{\bfv}(t)|_{{\bfH}_\sigma}. 
    \end{eqnarray*}
Next, we see from Lemma 3.2 that
for any $z\in H$
  \begin{eqnarray*}
   |f(\rho_\mu*w(t))z|_{V^*_0} 
    &=& \sup_{v \in V_0,|v|_{V_0} \le 1} 
     \langle f(\rho_\mu*w(t))z, v \rangle_0 
   = \sup_{v \in V_0,|v|_{V_0} \le 1} \langle z, 
     f(\rho_\mu*w(t))v
   \rangle_0 \\
 &\le & |z|_{V^*_0}
\sup_{v \in V_0,|v|_{V_0} \le 1}|f(\rho_\mu*w(t))v|_{V_0}
 \le M^\mu_1|z|_{V^*_0}.
  \end{eqnarray*}
Therefore,
$$ 
|f(\rho_\mu*w(t))z -bz|_{V^*_0}
    \le |f(\rho_\mu*w(t))z)|_{V^*_0} 
+ b|z|_{V^*_0} \le M^\mu_3 |z|_{V^*_0}. 
 $$
Thus \eqref{eq:3.8} is obtained.  \hfill $\Box$ \vspace{0.5cm}

\noindent
{\bf Proof of Proposition 3.1.} We shall prove the proposition in three steps.
\\[0.1cm]
(Step 1) Assume that ${\bfv} \in C([0,T]; {\bfV}_\sigma)$. By virtue of 
Lemma 3.3, our perturbation term
 $$h(t, z):= f(\rho_\mu*w(t))z -bz 
-{\bfv}(t)\cdot \nabla[\rho_\mu*(\gamma_\mu z)]
  $$ 
is Lipschitz continuous in $z\in V^*_0$ and continuous in $t$, so that it 
satisfies the condition $(h4)$ in Appendix III.
The other conditions $(h1)-(h3)$ are easily checked.
Therefore, the existence-uniqueness of a (strong) solution $u$ of 
$(B;w,{\bfv};u_0)_\mu$  
is a direct consequence of Proposition III; actually it admits one and
only one solution $u$ such that
$u \in W^{1,2}(0,T;V^*_0)$ and $t \to \varphi(u(t))=
|\hat \beta(u(t))|_{L^1(\Omega)}$ is absolutely continuous on $[0,T]$.
Since $0\le u \le u^*$ a.e.~on $Q$, these regularities imply 
$u \in C_w([0,T];H)$, where $ C_w([0,T];H)$ stands for the space of all
weakly continuous functions from $[0,T]$ into $H$.

Next, we show the uniform estimate \eqref{eq:3.6}. 
We observe that, by \eqref{eq:2.2}, and as  $|\rho_\mu|\le 1$,
 \begin{eqnarray*} 
  |f(\rho_\mu*w(t))u(t) -bu(t)|_{V^*_0}
 &\le& c_0|f(\rho_\mu*w(t))u(t) -bu(t)|_H\\
&\le& c_0\left\{(\max_{0\le r \le 1}f(r))u^*
  +bu^*\right\}|\Omega|^{\frac 12},
 \end{eqnarray*}
 and in the same way with the Remark 2.2 we obtain
 \begin{eqnarray*} 
  |{\bfv}(t)\cdot\nabla[\rho_\mu*(\gamma_\mu
   u(t))]|_{V^*_0}
&=&|{\rm div}[\rho_\mu*(\gamma_\mu u(t))
{\bfv}(t)]|_{V^*_0}\\
&=&\sup_{v\in V_0,|v|_{V_0}\le 1} \left\{-\int_\Omega
 \rho_\mu*(\gamma_\mu u(t)){\bfv}(t)\cdot \nabla v dx
    \right\}\\
&\le& u^*|{\bfv}(t)|_{{\bfH}_\sigma}
  \sup_{v \in V_0,|v|_{V_0}\le 1}|v|_{V_0}
  = u^*|{\bfv}(t)|_{{\bfH}_\sigma}.
\end{eqnarray*}
These inequalities imply that the perturbation
term $h(t,u(t))$ satisfies
  $$|h(\cdot,u)|_{L^2(0,T;V^*_0)} \le 
    M_4(1+|{\bfv}|_{L^2(0,T;{\bfH}_\sigma)})$$
for a positive constant $M_4$ independent of 
$\mu \in (0,1],~{\bfv}$ and $u$. 
Accordingly, from Appendix I, Proposition I(3), 
it follows that \eqref{eq:3.6}
holds for a non-negative increasing function 
$B_0(\cdot,\cdot)$.\vspace{0.2cm}

\noindent
(Step 2) In the general case of ${\bfv} \in L^2(0,T;{\bfV}_\sigma)$, we
choose a sequence $\{{\bfv}_n\}$ in 
$C([0,T]; {\bfV}_\sigma)$ such that 
${\bfv}_n \to {\bfv}$ in 
$L^2(0,T;{\bfV}_\sigma)$ (as $ n\to \infty$). According to the result of
(Step 1), $(B;w,{\bfv}_n;u_0)_\mu $ admits a unique solution 
$u_n$ which
enjoys the uniform estimate \eqref{eq:3.6}. Therefore we can choose a subsequence
$\{u_{n_k}\}$ from $\{u_n\}$ and a function $ u \in W^{1,2}(0,T;V^*_0)$ with
$\sup_{t \in[0,T]}|\hat\beta( u(t))|_{L^1(\Omega)}<\infty$ such that
  $$ u_{n_k} \to  u~{\rm in ~}C([0,T];V^*_0)
   ~{\rm and~weakly~in~}W^{1,2}(0,T;V^*_0),~\sup_{k\geq 1,~t \in [0,T]}
     |\hat\beta(u_{n_k}(t))|_{L^1(\Omega)} < \infty.  $$
Now it follows from Lemma 3.1 and Lemma 3.3 that 
$$f(\rho_\mu*w)u_{n_k} -b u_{n_k} 
-{\bfv}_{n_k} \cdot \nabla[\rho_\mu*
(\gamma_\mu u_{n_k})]
\to  f(\rho_\mu*w) u -b u -{\bfv}\cdot 
\nabla [\rho_\mu*(\gamma_\mu u)]~
 ~{\rm in~} L^2(0,T:V^*_0).$$
As a consequence, by Proposition II in the appendix, $u_{n_k}$ converges
in $C([0,T];V^*_0)$ to the solution of $(B;w,{\bfv},u_0)_\mu $. Clearly 
this solution coincides with $u$.\vspace{0.2cm}

\noindent
(Step 3) We now show uniqueness of solution. Let $u$ and $\bar u$ be
two solutions of $(B;w,{\bfv};u_0)_\mu$. Then it follows from the appendix,
Proposition I, (2), and from Lemma 3.3, that
 \begin{multline*}
 \frac 12 |u(t)-\bar u(t)|^2_{V^*_0}
     \ \le \ -\int_0^t (f(\rho_\mu*w)(u-\bar u)
    -b(u-\bar u) -{\bfv}\cdot 
   \nabla[\rho_\mu*(\gamma_\mu(u-\bar u))], 
          u-\bar u)_*d \tau\\
    \ \le \ (M^\mu_2|{\bfv}
    |_{L^\infty(0,T;{\bfH}_\sigma)}+M^\mu_3)
          \int_0^t|u-\bar u|^2_{V^*_0}d\tau.
  \end{multline*}
Therefore, by the Gronwall inequality, we have $u=\bar u$ on $[0,T]$.
\hfill  $\Box$ \vspace{0.5cm}

\noindent
{\bf Proposition 3.2.} {\it Assume  \eqref{eq:1.1} and
let $u_0\in H$ be such that $\hat{\beta}(u_0)\in L^1(\Omega)$.
Take any $\mu\in[0,1]$ and let
$\{\mu_n\}$ be a non--increasing sequence in $(0,1]$ 
such that $\mu_n\downarrow \mu$ (as $n\to\infty$).
Let $\{{\bfv}_n\}$ and $\{w_n\}$ be
sequences such that
   \begin{equation}\label{eq:3.9} 
   \left \{
      \begin{array}{l}
     \{{\bfv}_n\}~{\it is~bounded~in~}L^\infty(0,T;{\bfH}_\sigma),
   ~{\bfv}_n \to {\bfv}~{\it weakly~in}~L^2(0,T;{\bfV}_\sigma),\\[0.2cm]
   0\le w_n \le 1~ {\it a.e.~on~}Q,
   ~  w_n \to w~~{\it in ~}L^2(Q)~{\it and~weakly~in~}W^{1,2}(0,T;V^*), 
      \end{array} \right. 
   \end{equation}
(as $n \to \infty$).
Then $u_n$, the solution of 
$(B;w_n,{\bfv}_n,u_0)_{\mu_n} $, converges to the solution $u$ 
of $(B;w,{\bfv}, u_0)_\mu$ 
in the sense that
 \begin{equation}\label{eq:3.10} 
 u_n \to u~{\it in ~}C([0,T];V^*_0)\cap
 L^2(Q)~{\it and~weakly~in~}W^{1,2}(0,T;V^*_0),
 \end{equation}  
and
 \begin{equation}\label{eq:3.11} 
  \int_0^T \varphi(u_n(t))dt \to \int_0^T \varphi(u(t))dt. 
 \end{equation}
}
\noindent
\hspace*{-0.1cm}{\bf Proof.} We give the proof only in the case 
$\mu = 0$ (see Remark 3.1), the others being similar.
On account of the uniform estimate \eqref{eq:3.6}, $\{u_n\}$ is bounded in
$W^{1,2}(0,T;V^*_0)$ and $0\le u_n \le u^*$ a.e.~on $Q$. Therefore
there is a subsequence of $\{u_n\}$, 
that we still denote  by $\{u_n\}$,  such that
$u_{n} \to u$ in $C([0,T];V^*_0)$ (as $n\to \infty$) 
for a certain function $u$ satisfying the estimate \eqref{eq:3.6}. Now, 
put $g_n(t):= f(\rho_{\mu_n}*w_n(t))\,u_n(t)-b\,u_n(t)-{\bfv}_n(t)\cdot 
\nabla [\rho_{\mu_n}*(\gamma_{\mu_n} u_n(t))]$ and
$g(t):= f(w(t))u(t)-bu(t)-{\bfv}(t)\cdot \nabla u(t)$.
Since $\{g_n\}$ is bounded in $L^2(0,T;V^*_0)$, 
it follows from Lemma 3.1 that $\{u_n\}$ 
is relatively compact in $L^2(Q)$, and hence
converges to~$u$ in $L^2(Q)$. This shows that $\gamma_{\mu_n}u_n\to u$
in~$L^2(Q)$ as well as 
$\rho_{\mu_n}* (\gamma_{\mu_n}u_n)\to u$ in~$L^2(Q)$. 
Besides,
$g_{n}\to g$ weakly in $L^2(0,T;V^*_0)$, which is seen as follows.
Observe that
\begin{eqnarray*}
 g_{n}-g
&=& (f(\rho_{\mu_n}*w_n)-f(w))u_n
+ f(w)(u_n-u)
     -b\,(u_n-u)\\
& &-{\bfv} \cdot \nabla [\rho_{\mu_n}*(\gamma_{\mu_n}u_n)-u] 
-(\bfv_n-\bfv)\cdot\nabla[\rho_{\mu_n}*(\gamma_{\mu_n} u_n)].
\end{eqnarray*}
From the assumption \eqref{eq:3.9} with \eqref{eq:2.6}--\eqref{eq:2.9}
it follows that the first four terms at  the right hand side
converge to $0$ in
$C([0,T]; V^*_0)$, and the last one converges weakly  to $0$ in $L^2(0,T;V^*_0)$. Therefore, the limit $u$ is
a unique solution of $(B;w, {\bfv},u_0)$, and  \eqref{eq:3.10} and \eqref{eq:3.11} hold 
by Proposition~II in Appendix~II. \hfill $\Box$

\section{Nutrient transport equation and its approximation}\label{sec:nutrient}

Given functions $u \in C_w([0,T];H)$ with 
$0\le u \le u^*$ a.e.~on $Q$
and ${\bfv} \in L^2(0,T;{\bfV}_\sigma)$, our nutrient transport equation is 
treated in the form:
   $$(N;u,{\bfv},w_0) ~~\left \{
     \begin{array}{l}
      w'(t) +\partial_{V^*} \Phi^t(u;w(t))
          +{\bfv}(t)\cdot \nabla w(t)
          = -f(w(t))u(t)\\[0.12cm]
~~~~~~~~~~~~~~~~~~~~~~~~~~~~~~~~~~~~~~~~~~~~~~~~~~~%
{\rm in~}V^*~{\rm for~a.e.~}t \in [0,T],\\
    w(0)=w_0,
     \end{array} \right. $$
where the initial datum
$w_0$ is prescribed in $H$, satisfying
$0\le w_0\le 1$ a.e.~on $\Omega$, 
$f(\cdot)$ satisfies \eqref{eq:1.3}, 
and $\Phi^t(u; \cdot)$ is a non-negative, continuous
and convex function on $V$ defined by
  $$ \Phi^t(u; w):=
       \frac 12 \int_\Omega d(u(x,t))
         |\nabla w(x)|^2 dx,~~~\forall w \in V,
  $$
with the function $d(\cdot)$ satisfying \eqref{eq:1.2}; $\partial_{V^*}\Phi^t(u; \cdot)$
is the subdifferential of $\Phi^t(u;\cdot)$ from $V=D(\partial_{V^*}
\Phi^t(u;\cdot))$ into $V^*$. We see that $\partial_{V^*}\Phi^t(u; \cdot)$
is singlevalued, linear and maximal monotone from $V$ into $V^*$, satisfying
 $$
     \langle \partial_{V^*}\Phi^t(u; w),z \rangle 
     =\int_\Omega d(u(x,t))\nabla w(x)\cdot \nabla z(x)dx,
   ~~\forall w, ~z \in V,~\forall t \in [0,T].$$

\noindent
{\bf Definition 4.1.} {\it Let $u \in C_w([0,T];H)$
with $0\le u \le u^*$ a.e.~on $Q$
and ${\bfv} \in L^2(0,T;{\bfV}_\sigma)$.
Then, for 
any $w_0 \in H$ with $0\le w_0\le 1$
a.e.~on $\Omega$, a function $w: [0,T] \to V$ is called a solution
of $(N; u,{\bfv}, w_0)$, if $w \in L^2(0,T;V)\cap L^\infty(Q)$,  
$w'\in L^2(0,T;V^*)$, $w(0)=w_0$ and
 \begin{equation}\label{eq:4.1} 
 w'(t) +\partial_{V^*} \Phi^t(u;w(t))
     +{\bfv}(t)\cdot \nabla w(t)= -f(w(t))u(t)~{\rm in ~}V^*~
  {\rm for~a.e.~}t\in [0,T].
\end{equation}
}

\noindent
{\bf Remark 4.1.} We shall construct
a solution $w$ such that $0\le w \le 1$ a.e.~on $Q$.
\vspace{0.25cm}

\noindent
{\bf Remark 4.2.}
If $w \in L^\infty(Q)$ or ${\bfv} \in L^2(0,T;\bfH_\sigma)\cap L^\infty(Q)^3$ , 
we have 
$\nabla w\cdot\bfv={\rm div}(w{\bfv})\in L^2(0,T;V^*)$, cf.~Remark~2.2.
Indeed, assume $w \in\L^\infty(\Omega)$, then, for all $z \in L^2(0,T;V)$:
$$
\int_0^T \langle {\rm div}(w{\bfv}), z\rangle dt 
=-\int_0^T \int_\Omega w{\bfv} \cdot \nabla z dxdt 
\le |w|_{L^\infty(Q)}\,|{\bfv}|_{L^2(0,T;{\bfH}_\sigma)}\,|z|_{L^2(0,T;V)}.
$$
The other case is analogous.
\vspace{0.25cm}

\noindent
{\bf Remark 4.3.}
If ${\bfv} \in L^2(0,T;\bfH_\sigma)\cap L^\infty(Q)^3$, the linear operator
$w \to {\bfv}(t)\cdot \nabla w $
is continuous from $V$ into $V^*$ 
and  maximal monotone. 
Indeed, by Remark~2.2, 
   \begin{equation}\label{eq:4.2} 
   \int_\Omega \left({\bfv}(x,t)\cdot \nabla w(x)\right) w(x) dx
     =\frac12\int_\Omega {\rm div}\left({\bfv}(x,t) w(x)^2\right) dx=0 
 \end{equation}
for all $w \in V$ and $t \in [0,T]$. Therefore the sum
$w \to \partial_{V^*}\Phi^t(u;w)+{\bfv}(t)\cdot \nabla w$ is linear,
continuous, maximal monotone and coercive from $V$ into $V^*$.\vspace{0.5cm}

\noindent
We recall  the general theory on evolution equations with monotone
operators in Banach spaces (cf. [3; Chapter 4]) for the solvability
of $(N; u,{\bfv},w_0)$. On account of Remark~4.3, this gives 
the following lemma. \vspace{0.5cm}

\noindent
{\bf Lemma 4.1.} {\it Assume that $u \in C_w([0,T];H)$ with $0\le u \le u^*$,
${\bfv} \in L^2(0,T;{\bfV}_\sigma)\cap L^\infty(Q)$ 
and \eqref{eq:1.2} is satisfied. 
Then we have:
\begin{description}
\item{(1)} For any $f^* \in L^2(0,T;V^*)$ and $w_0 \in H$
the Cauchy problem
 \begin{equation}\label{eq:4.3} 
  \left\{
    \begin{array}{l}
  w'(t)+\partial_{V^*}\Phi^t(u;w(t))+{\bfv}(t)\cdot \nabla w(t)
  =f^*(t)~~{\rm in~}V^*~{\rm for~a.e.~}t \in [0,T],
  \\[0.2cm]
   w(0)=w_0,
     \end{array} \right.
     \end{equation}
admits one and only one solution $w$ such that
$ w\in L^2(0,T;V)$ and $w' \in L^2(0,T;V^*)$. 
\item{(2)}
Let $w_i$ be the solution of \eqref{eq:4.3}
with $w_0=w_{i0}\in H$ and $f^*=f^*_i \in 
L^2(0,T;V^*)$ for $i=1,2$. Then, for all $t\in [0,T]$:
 \begin{equation}\label{eq:4.4}
 \frac 12 |w_1(t)-w_2(t)|^2_H +
    c_d \int_0^t|\nabla(w_1-w_2)|^2dxd\tau
  \ \le \  \frac 12 |w_{10}-w_{20}|^2_H
    +\int_0^t \langle f^*_1-f^*_2, 
      w_1-w_2 \rangle\, d\tau,
    \end{equation} 
\end{description} 
}

\noindent
We prove now the existence-uniqueness result
for $(N;u,{\bfv},w_0)$.
 \vspace{0.25cm}

\noindent
{\bf Proposition 4.1.} {\it Assume that
$u\in C_w([0,T];H)$ with $0\le u\le u^*$ a.e.~on~$Q$, 
\eqref{eq:1.2}, \eqref{eq:1.3} are satisfied, 
${\bfv}\in L^2(0,T;{\bfV}_\sigma)$,
and $w_0 \in H$ with $0\le w_0 \le 1$ a.e.~on~$\Omega$. 
Then the problem
$(N;u,{\bfv},w_0)$ admits one and only one solution $w$. 
This solution satisfies 
  $$ 0 \le w \le 1~~{\it a.e.~on~}Q,$$
and
 \begin{equation}\label{eq:4.5}  
 |w(t)|^2_H+2c_d\int_0^t |\nabla w|^2_Hd\tau
   \le e^{2u^*L(f)T} |w_0|^2_H,~~\forall
    t \in [0,T]. 
 \end{equation}
 }
\noindent
{\bf Proof.} We prove the proposition in three steps.\vspace{0.2cm}

\noindent
(Step 1) Assume first that 
${\bfv} \in L^\infty(Q)^3$. 
We are going to construct the solution $w$ of $(N;u,{\bfv},w_0)$ by the contraction mapping principle. 
Let $T_1\in (0, T]$ be a time such that $2u^*L(f)T_1 <1$ and, using Lemma~4.1, 
define a mapping ${\cal N}:
C([0,T_1];H) \to C([0,T_1];H)$, 
which assigns to
each $\bar w \in C([0,T_1];H)$ the solution $w$ of \eqref{eq:4.3}
on $[0,T_1]$ with $f^*=f(\bar w)u$, namely 
$w:={\cal N}\bar w$. Then, for any $\bar w_i \in C([0,T_1];H),~i=1,2$, 
we observe from 
\eqref{eq:4.4}
 that 
\begin{multline*}
 \frac 12 |w_1(t)-w_2(t)|^2_H +
   c_d \int_0^t\int_\Omega|\nabla(w_1-w_2)|^2dxd\tau \\ 
 \le u^*\int_0^t |f(\bar w_1)-f(\bar w_2)|_H 
      |w_1-w_2|_H d\tau
 \le u^*L(f)\int_0^t|\bar w_1-\bar w_2|_H 
    |w_1-w_2|_H d\tau,
 \end{multline*}
for all $ t\in [0,T_1]$, so that 
 $$ |w_1-w_2|_{C([0,T_1];H)} \le 2 u^*L(f)T_1
                 |\bar w_1-\bar w_2|_{C([0,T_1];H)}.
$$                 
This shows that ${\cal N}$ is strictly contractive in $C([0,T_1];H)$
and it has a unique fixed point $w$ in $C([0,T_1];H)$, namely
$w={\cal N}w$, which is a unique solution of  \eqref{eq:4.3} 
on the time interval
$[0,T_1]$. It is a routine work to construct a unique solution $w$ of 
$(N;u,{\bfv};w_0)$
on the whole interval $[0,T]$ by a finite number of time-steps. 
\vspace{0.2cm}

\noindent
(Step 2) Still assume that
${\bfv} \in L^\infty(Q)^3$, and recall that $0\le w_0 \le 1$ a.e.~on $\Omega$. 
Then we show that the solution $w$ of \eqref{eq:4.3} constructed in 
Step~1 satisfies 
$ 0\le w \le 1$ a.e.~on~$Q$.
To do so, multiply \eqref{eq:4.1} by $-w^{-}$ ($=$ the negative part of $w$) and 
integrate the both sides in time to get by \eqref{eq:1.3} 
 $$ \frac 12 |w^{-}(t)|^2_H + c_d\int_0^t\int_\Omega
  |\nabla w^{-}|^2 dxd\tau \le u^*L(f)\int_0^t |w^{-}(\tau)|^2_Hd\tau,
   ~~\forall t \in [0,T].$$
Applying the Gronwall's lemma to this inequality, we obtain that 
$|w^{-}(t)|_H=0$ for all $t \in [0,T]$, namely $w \geq 0$ a.e.~on $Q$.
Similarly, by multiplying \eqref{eq:4.1} by $(w-1)^+$ ($=$ the positive part of $w-1$),
and integrating the both sides in time, 
we conclude that $|(w-1)^+|_H=0$, namely $ w \le 1$ a.e.~on $Q$. 
Thus $0\le w \le 1$ a.e.~on $Q$, and $w$ is the solution of $(N;u,{\bfv},w_0)$
in the sense of Definition~4.1. \vspace{0.2cm}

\noindent
(Step 3) For general ${\bfv}$, we
approximate ${\bfv}\in L^2(0,T; {\bfV}_\sigma)$ 
by a sequence $\{{\bfv}_n\}$ from  $L^2(0,T; {\bfV}_\sigma)\cap L^\infty(Q)^3$ such that ${\bfv}_n \to {\bfv}$
in $L^2(0,T; {\bfV}_\sigma)$ (as $n \to \infty$). 
By virtue of Steps 1 and 2, for each $n$, the problem 
 \begin{equation}\label{eq:4.6}  
 \begin{array}{l}\displaystyle
   w'_n(t) + \partial_{V^*}\Phi^t(u; w_n(t))
     +{\bfv}_n(t)\cdot \nabla w_n(t)= -f(w_n)u(t)~~{\rm in ~}V^*, 
\\[0.12cm]
w_n(0)=w_0, 
 \end{array}
 \end{equation}
 has a unique solution $w_n$ such that
$w_n \in L^2(0,T;V)$, $w'_n \in L^2(0,T; V^*)$ and $0\le w_n \le 1$ a.e.~on 
$Q$. Multiplying \eqref{eq:4.6} by $w_n$, we obtain 
by \eqref{eq:4.2} that
\begin{equation}\label{eq:4.7}   
   \frac 12 |w_n(t)|^2_H + c_d \int_0^t |\nabla w_n(\tau)|^2_H d\tau
     \le \frac 12 |w_0|^2_H +u^*L(f)\int_0^t |w_n|^2_H d\tau,~~\forall
    t \in [0,T], %
\end{equation}
which implies, with the Gronwall inequality, 
that $\{w_n\}$ is bounded in $L^2(0,T;V)$.
We also infer from ^^ ^^ $0\le w_n \le 1$" and Remark~4.2 that 
${\bfv}_n\cdot \nabla w_n={\rm div}(w_n {\bfv}_n)$ is bounded in 
$L^2(0,T;V^*)$. Consequently, by \eqref{eq:4.6}, $\{w'_n\}$ is bounded in $L^2(0,T;V^*)$.
Therefore there exist a subsequence $\{w_{n_k}\}$ of $\{w_n\}$ 
and a function $w \in L^2(0,T;V)$ with 
$0\le w\le 1$ a.e.~on~$Q$, such that $w_{n_k} \to w$ weakly in $L^2(0,T;V)$. Furthermore, on 
account of the Aubin's compactness theorem [2], 
we have $w_{n_k} \to w$ in $L^2(Q)$. Now it is easy to see,
by letting $k \to \infty$ in \eqref{eq:4.6} with $n=n_k$, that the limit $w$
satisfies \eqref{eq:4.1} and the same type of energy inequality 
as \eqref{eq:4.7} holds for $w$. 
We easily get the estimate \eqref{eq:4.5} from it. 
Uniqueness of solution and \eqref{eq:4.5} 
are obtained by the Gronwall inequality.
\hfill $\Box$
\vspace{0.5cm}

\noindent
{\bf Proposition 4.2.} {\it Assume that \eqref{eq:1.2}
and \eqref{eq:1.3} hold, $w_0 \in H$ with
$0\le w_0 \le 1$ a.e.~on~$\Omega$.
Let $\{{\bfv}_n\}$ and $\{u_n\}$ be sequences 
such that $0 \le u_n \le u^*$ a.e.~on $Q$ for
all $n$, and 
 \begin{equation}\label{eq:4.8}  
  {\bfv}_n \to {\bfv}~{\it~weakly~in~} L^2(0,T;{\bfV}_\sigma),~~
     u_n \to u~{\it in~}C([0,T];V^*_0)\cap L^2(Q). 
 \end{equation}
Then, the solution $w_n$ of $(N;u_n,{\bfv}_n,w_0)$ converges to the solution 
$w$ of $(N;u,{\bfv},w_0)$ in the sense that
 \begin{equation}\label{eq:4.9} 
 w_n \to w~{\it in~}L^2(Q)~{\it and~weakly~in~}L^2(0,T;V),~
 w'_n \to w'~{\it weakly~in~}L^2(0,T;V^*). 
\end{equation}
} 
\noindent
\hspace*{-0.1cm}{\bf Proof.} From the uniform estimate \eqref{eq:4.5}
we observe that $\{w_n\}$ is bounded in $L^2(0,T;V)$ with
$0 \le w_n\le 1$ a.e.~on $Q$, 
so that $w'_n= -\partial_{V^*}\Phi^t(u_n;w_n)- 
{\rm div}({\bfv}_n w_n)- f(w_n)u_n$
is bounded in $L^2(0,T;V^*)$. It follows  from the Aubin's
compactness theorem [2] that $\{w_n\}$ is relatively compact in $L^2(Q)$.
Therefore, there are a subsequence $\{w_{n_k}\}$ of $\{w_n\}$
and a function $\bar w$ so that
$w_{n_k} \to  \bar w~~{\rm in~}L^2(Q)$
and weakly in $L^2(0,T:V)$ as well as
$w'_{n_k} \to \bar w'$ weakly in 
$L^2(0,T;V^*)$.
By these convergences and \eqref{eq:4.8} we see that 
$\partial_{V^*}\Phi^t(u_{n_k};w_{n_k}) \to \partial_{V^*}\Phi^t(u;\bar w)$
weakly in $L^2(0,T;V^*)$ and
$
{ -{\rm div}(w_{n_k}{\bfv}_{n_k}) - f(w_{n_k})u_{n_k}} 
\to
{ -{\rm div}(\bar w{\bfv}) - f(\bar w)u }
$ 
weakly in 
$L^2(0,T;V^*)$ (as $k\to \infty$).
Hence, by Remark 4.2,  the limit $\bar w$ is a solution of $(N;u,{\bfv},w_0)$. 
By uniqueness we have $\bar w= w$, which
implies that convergences \eqref{eq:4.9} hold without extracting any subsequence
from $\{w_n\}$. \hfill $\Box$
\vspace{0.3cm}

\noindent
As a regular approximation for $(N; u,{\bfv},w_0)$, we employ problem
$(N; \rho_\mu*u,{\bfv},w_0)$, which is denoted by $(N; u,{\bfv},w_0)_\mu$
for any small parameter $\mu \in (0,1)$,
namely
$$ (N; u,{\bfv},w_0)_\mu~~\left\{
  \begin{array}{l}
    w'(t)+\partial_{V^*}\Phi^t(\rho_\mu*u; w(t))+{\bfv}\cdot \nabla w(t)
     =-f(w(t))\rho_\mu*u\\
     ~~~~~~~~~~~~~~~~~~~~~~~~{\rm in~}V^*~{\rm for~a.e.~} t \in [0,T],\\
    w(0)=w_0.
   \end{array} \right. $$
It is clear that Propositions 4.1,  4.2 are valid for this approximate
problem by replacing $u$ by $\rho_\mu*u$.

\section{Variational inequality of the Navier-Stokes type and its approximation}\label{sec:hydro}
As was mentioned in the introduction, the biomass formation mechanism, 
together with the  nutrient transport and consumption
takes place in a fluid. 
At the same time, the forming biomass becomes an obstacle for the flow. We model it  by making use of a variational
inequality of Navier-Stokes type. 

Let $p_0: (0,u^*] \to {\bf R}$ be the same function as in (i) in the introduction, satisfying \eqref{eq:1.1} (see Fig.1), 
and let $u$ be a given function in 
$C_w([0,T];H)$ with 
$0\le u \le u^*$ a.e.~on $Q$. Then, with the function 
$u^\varepsilon:=\rho_\varepsilon*u$ for a fixed small positive parameter 
$\varepsilon \in (0,1)$, the strong
formulation of our variational inequality of the Navier-Stokes type is of the
following form:
$$(H;u, {\bfv}_0, \bfg)^\varepsilon ~~~ \left\{
    \begin{array}{l}
 \displaystyle{|{\bfv}(x,t)|\le p_0(u^\varepsilon(x,t))
  ~~{\rm for~}(x,t) \in Q;}\\[0.3cm]
   \displaystyle{\langle {\bfv}'(t), {\bfv}(t) - {\bfz}\rangle_\sigma
  +\nu\int_\Omega \nabla {\bfv}(x,t)\cdot \nabla({\bfv}(x,t)-{\bfz(x)})dx}
   \\[0.3cm]
   \displaystyle{~~~~~~~~+\int_\Omega({\bfv}(x,t)\cdot \nabla)
     {\bfv}(x,t)\cdot({\bfv}(x,t)-{\bfz}(x))dx 
     \le (\bfg(t), \bfv(t)-\bfz)_\sigma} 
     \\[0.5cm]
  \displaystyle{~~~~\forall {\bfz} \in {\bfV}_\sigma~{\rm with~}
  |{\bfz}(x)|\le p_0(u^\varepsilon(x,t)) ~{\rm for~}x \in \Omega,
   ~t \in[0,T],}\\[0.3cm]
  \displaystyle{ {\bfv}(x,0)={\bfv}_0(x)~~{\rm for ~}x \in \Omega,}
      \end{array}
     \right. $$
where $\nu$ is positive constant (viscosity), 
${\bfv}_0$ a prescribed initial datum for ${\bfv}$ 
and $\bfg\in L^2(0,T; \bfH_\sigma)$ a prescribed external force.  
\vspace{0.25cm}

The existence-uniqueness of a strong solution 
to $(H;u,{\bfv}_0, \bfg)^\varepsilon$ 
is of course an open question just as the usual 3D Navier-Stokes equations.
Therefore, we shall construct a weak solution 
of $(H;u,{\bfv}_0, \bfg)^\varepsilon$ in the variational sense.
\vspace{0.5cm}

\noindent
{\bf Definition 5.1.} {\it Let $u \in C_w([0,T];H)$
with $0\le u\le u^*$ a.e.~on $Q$,~
$u^\varepsilon:=\rho_\varepsilon*u$ and
${\boldsymbol{\mathcal K}}(u^\varepsilon)$ be the class of
test functions defined by
 \begin{equation*}\label{eq:5.1}
 {\boldsymbol{\mathcal K}}(u^\varepsilon):=\left\{{\boldsymbol{\mathcal \eta}} \in C^1([0,T];{\bfW}^{1,4}_{0,\sigma}
    (\Omega))~\left|~
       \begin{array}{l}
     \text{supp}\,({\boldsymbol{\mathcal \eta}})\subset 
     \hat Q(u^\varepsilon < \delta_0),\\[0.2cm]
     |{\boldsymbol{\mathcal \eta}}| \le p_0(u^\varepsilon)~\text{on~}Q
       \end{array}
        \right. \right\} 
  \end{equation*}
where ${\bfW}^{1,4}_{0,\sigma}(\Omega)$ is the closure of 
${\bf{\mathcal D}}_\sigma(\Omega)$ in $W^{1,4}_0(\Omega)^3$ and
$\hat Q(u^\eps <\delta_0):=\{(x,t) \in \Omega\times [0,T]~|~u^\eps(x,t) 
<\delta_0\}$.
Then, for a given initial datum ${\bfv}_0$ 
and  $\bfg\in L^2(0,T; \bfH_\sigma)$, 
a function ${\bfv}:[0,T]\to {\bfH}_\sigma$
is called a weak solution of 
$(H;u,{\bfv}_0,\bfg)^\varepsilon$, if the following
conditions hold:
\begin{description}
\item{(1)} ${\bfv} \in L^2(0,T;{\bfV}_\sigma)$
{and $\sup_{t \in [0,T]}|{\bfv}(t)|_{{\bfH}_\sigma}<\infty$};
\item{(2)} the function
$t \to ({\bfv}(t),
{\boldsymbol{\mathcal \eta}}(t))_\sigma$ is of bounded variation on $[0,T]$ for
any ${\boldsymbol{\mathcal \eta}} \in{\boldsymbol{\mathcal K}}(u^\varepsilon)$;
\item{(3)} ${\bfv}$ satisfies: $ \qquad {\bfv}(0)={\bfv}_0,~~ |{\bfv}(x,t)|\le p_0(u^\varepsilon(x,t))~
 {\rm a.e.~}x \in \Omega,~\forall t \in [0,T],$
\begin{equation}\label{eq:5.2}
\begin{split}
\int_0^t ({\boldsymbol{\mathcal \eta}}'(\tau), {\bfv}(\tau)
   -{\boldsymbol{\mathcal \eta}}(\tau))_\sigma \,d\tau
  +\nu \int_0^t\int_\Omega \nabla {\bfv}(x,\tau)\cdot\nabla({\bfv}(x,\tau)
  -{\boldsymbol{\mathcal \eta}}(x,\tau))\,dx d\tau 
  \\
   +\int_0^t\int_\Omega ({\bfv}(x,\tau)\cdot \nabla){\bfv}(x,\tau)\cdot
     ({\bfv}(x,\tau)-{\boldsymbol{\mathcal \eta}}(x,\tau))\, dx  d\tau
   +\frac 12 |{\boldsymbol{\mathcal \eta}}(t)-{\bfv}(t)|^2_{{\bfH}_\sigma}\\
   \le \int_0^t (\bfg(\tau), \bfv(\tau) - \eta(\tau))_\sigma\, d\tau 
   +
   \frac 12 |{\boldsymbol{\mathcal \eta}}(0)-{\bfv}_0|^2_{{\bfH}_\sigma}, ~~~~~~~\\
 \forall t \in [0,T],~\forall {\boldsymbol{\mathcal \eta}}
 \in{\boldsymbol{\mathcal K}}(u^\varepsilon).
\end{split}
 \end{equation}
\end{description}
 }
 \noindent
In the rest of this section, we propose an 
approximate problem 
$(H;u,{\bfv}_0, \bfg)^\varepsilon_\mu$ for 
$(H;u,{\bfv}_0, \bfg)^\varepsilon$.
We begin with the approximation $p_\mu(r)$ of $p_0(r)$ with a small positive
parameter $\mu$ (actually $\mu \in (0,\delta_0)\cap (0,1)$ with 
$\mu <p_0(\mu)$), See Fig.2:
\begin{equation}\label{eq:5.3}
	 p_\mu(r):=\left \{
        \begin{array}{ll} 
        p_0(\mu),~~&{\rm for~}r \in [0,\mu],\\[0.2cm]
        p_0(r),~~&{\rm for~}r \in (\mu, p_0^{-1}(\mu)],\\[0.2cm]
        \mu,~~&{\rm for~}r \in (p_0^{-1}(\mu), u^*].
        \end{array}
       \right. 
\end{equation}       
\begin{figure}[b]
\begin{center}
\vspace*{-.5cm}
\includegraphics[width=14cm, height=4.5cm]{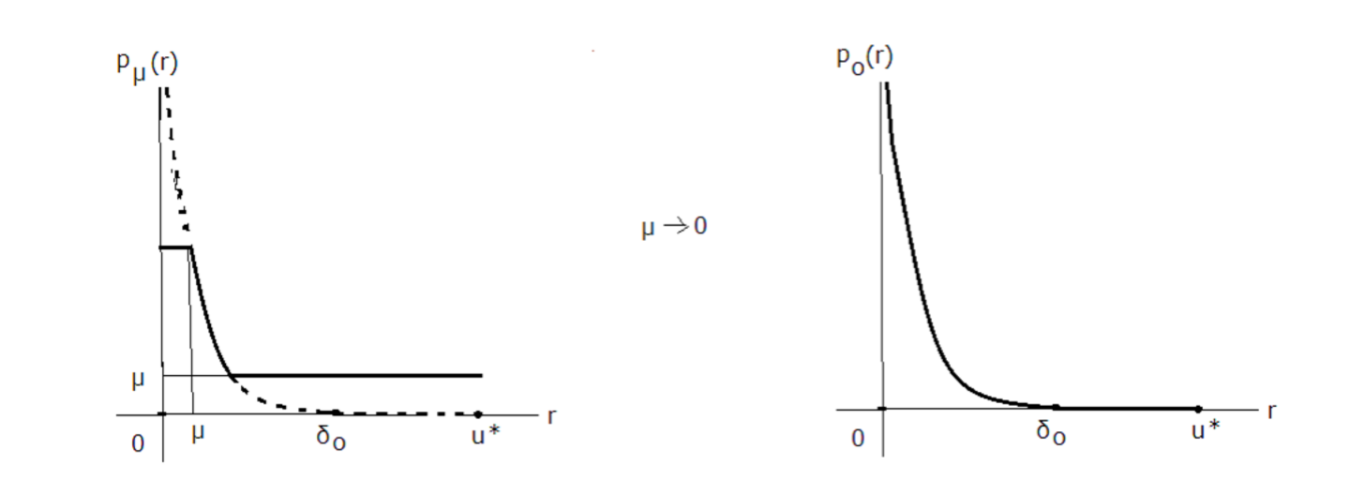}
\end{center}
\vspace*{-1cm}
\caption{Approximating the obstacle function $p_0$ by $p_\mu$ }
\end{figure}
Next, we approximate the obstacle function $p_0(u^\varepsilon)$ by
$p_\mu((\gamma_\mu u)^\varepsilon)$ with $(\gamma_\mu u)^\varepsilon:=
\rho_\varepsilon*(\gamma_\mu u)$, 
where $\gamma_\mu$ is given by \eqref{eq:3.4}. We put finally
  $$K_{\mu}((\gamma_\mu u)^\varepsilon;t)
    :=\{{\bfz} \in {\bfV}_\sigma~|~
  |{\bfz}(x)| \le p_{\mu}((\gamma_\mu u)^\varepsilon(x,t))~{\rm a.e.~}x 
    \in \Omega\},~~ \forall t\in [0,T]. $$
Now, consider the following approximate problem 
$(H;u,{\bfv}_0,\bfg)^{\varepsilon}_{\mu}$ for $(H;u,{\bfv}_0,\bfg)^\varepsilon$:
 $$ (H;u,{\bfv}_0, \bfg)^{\varepsilon}_{\mu}~\left \{
   \begin{array}{l}
   {\bfv}(t) \in K_{\mu}((\gamma_\mu u)^\varepsilon; t),~\forall t \in [0,T];
\\[0.3cm]
   \displaystyle{\langle {\bfv}'(t), {\bfv}(t)-{\bfz}\rangle_\sigma
     +\nu\int_\Omega 
   \nabla {\bfv}(t)\cdot \nabla({\bfv}(t)-{\bfz}) dx ~~~~~~~~~~~~~~~~~}
   \\[0.3cm]
    ~~~~~~~~~~~~~~~~~~~~~~~~  
    +\langle {\bf{\mathcal G}}({\bfv}(t),{\bfv}(t)),{\bfv}(t)-{\bfz}\rangle_\sigma
    \le (\bfg(t), \bfv(t) - \bfz)_\sigma\\[0.3cm]
    ~~~~~~~~~~~~~~~~~~~~~~~~~~~~~~~~~~~~~~~~~~~~~~\forall {\bfz}
     \in K_{\mu}((\gamma_\mu u)^\varepsilon; t),~{\rm a.e.~} t \in [0,T],
      \\[0.3cm]
     {\bfv}(0)={\bfv}_0;
   \end{array} \right. $$
 where ${\bf{\mathcal G}}$ is a nonlinear operator from 
 ${\bfV}_\sigma \times {\bfV}_\sigma \to {\bfV}^*_\sigma$
given by
 $$ \langle {\bf{\mathcal G}}({\bfv},{\bfw}),{\bfz}\rangle_\sigma
    :=\sum_{k,j=1}^3 \int_\Omega v^{(k)}\frac{\partial w^{(j)}}{\partial x_k}
     z^{(j)}dx $$
for all $\bfv:=(v^{(1)},v^{(2)},v^{(3)})$, 
$ {\bfw} :=(w^{(1)},w^{(2)},w^{(3)})$ 
and 
${\bfz}:=(z^{(1)},z^{(2)},z^{(3)})$ in ${\bfV}_\sigma \cap L^\infty(\Omega)^3$. 

\bigskip

\noindent
{\bf Remark 5.1} (a) By divergencee freeness of ${\bfv}\in\bfV_\sigma$, we have
$\langle {\bf{\mathcal G}}({\bfv},{\bfv}), {\bfv}\rangle_\sigma =0.$ \\[0.12cm]
(b) Also, 
${\bf{\mathcal G}}({\bfv},{\bfv}) \in {\bfH}_\sigma$ for
${\bfv} \in K_\mu((\gamma_\mu u)^\varepsilon;t)$.
\bigskip

In order to describe the above variational inequality as an evolution
inclusion of the subdifferential type we introduce 
time-dependent convex functions, $\psi_\mu^t((\gamma_\mu u)^\varepsilon; \cdot)$, on ${\bfH}_\sigma$,  of the following form:
 $$ \psi_\mu^t((\gamma_\mu u)^\varepsilon; {\bfz}):= \left\{
       \begin{array}{ll}
        \displaystyle{ \frac \nu 2 |{\bfz}|^2_{{\bfV}_\sigma},}
           &~~~{\rm if~}{\bfz} 
                    \in K_\mu((\gamma_\mu u)^\varepsilon; t),\\[0.5cm]
         \infty,&~~~{\rm otherwise}
        \end{array} \right. $$
and denote by 
$\partial\psi_\mu^t((\gamma_\mu u)^\varepsilon; \cdot)
= \partial_{\bfH_\sigma}\psi_\mu^t((\gamma_\mu u)^\varepsilon; \cdot)$ their
subdifferential
in ${\bfH}_\sigma$. We see that
${\bfv}^* \in \partial \psi^t_\mu((\gamma_\mu u)^\varepsilon; {\bfv})$ 
if and only if ${\bfv} \in K_\mu((\gamma_\mu u)^\varepsilon;t)$, 
${\bfv}^* \in {\bfH}_\sigma$ and
\begin{equation}\label{eq:5.sub}
({\bfv}^*, {\bfz} - {\bfv})_\sigma \le \nu \int_\Omega \nabla {\bfv}(x)
   \cdot \nabla({\bfz}(x)-{\bfv}(x))dx,~~\forall {\bfz} \in 
    K_\mu((\gamma_\mu u)^\varepsilon;t). 
\end{equation}
 
\noindent 
Now, we take ${\bfv}_0\in K_\mu((\gamma_\mu u)^\varepsilon;0)$ 
and consider the evolution inclusion:
\begin{equation}\label{eq:5.4}
 \left\{
    \begin{array}{l}
  \displaystyle{  
    {\bfv}'(t)+\partial \psi^t_\mu((\gamma_\mu u)^\varepsilon; {\bfv}(t)) 
  +{\bf{\mathcal G}}({\bfv}(t),{\bfv}(t)) \ni \bfg(t)~~{\rm in~}{\bfH}_\sigma
    ~{\rm for~a.e.~}t \in [0,T],}\\[0.12cm]
  \displaystyle{ {\bfv}(0)={\bfv}_0,}
    \end{array} \right. 
\end{equation}     
    By Remark~5.1(b), 
    \eqref{eq:5.4} makes sense as an inclusion in~${\bfH}_\sigma$. 
If 
${\bfv} \in W^{1,2}(0,T;{\bfH}_\sigma)$, then
\eqref{eq:5.4} is equivalent to $(H;u,{\bfv}_0, \bfg)^\varepsilon_\mu$ 
by \eqref{eq:5.sub}.  
A function ${\bfv}: [0,T] \to {\bfH}_\sigma$ is called a (strong) solution
to $(H;u,{\bfv}_0, \bfg)^\varepsilon_\mu$, if ${\bfv} \in W^{1,2}(0,T;{\bfH}_\sigma)
\cap C([0,T];{\bfV}_\sigma)$ and \eqref{eq:5.4} holds.
\vspace{0.5cm}

\noindent
{\bf Proposition 5.1.} {\it Let $\mu$ be any small positive number, and let 
$u$ be any function in $W^{1,2}(0,T;V^*_0)$ with $0\le u\le u^*$ a.e.~on $Q$
(hence $u \in C_w([0,T];H)$), let $\bfg\in L^2(0,T;\bfH_\sigma)$. 
Also, let ${\bfv}_0$ be any function in 
$K_{\mu}((\gamma_\mu u)^\varepsilon;0)$.
Then
$(H;u,{\bfv}_0,\bfg)^{\varepsilon}_{\mu}$ 
has one and only one solution ${\bfv}$, 
satsfying
\begin{equation}\label{eq:5.5}
 |{\bfv}(t)|^2_{{\bfH}_\sigma}+\nu \int_0^t |{\bfv}(\tau)|^2
   _{{\bfV}_\sigma}d\tau
   \le |{\bfv}_0|^2_{{\bfH}_\sigma}
   +\frac{L_P^2}{\nu}\int_0^T|\bfg(\tau)|^2_{\bfH_\sigma} d\tau,
\end{equation} 
where $L_p$ is the Poincar\'e constant, 
i.e.~$|\bfz|_{\bfH_\sigma}\le L_P|\bfz|_{\bfV_\sigma}$ 
for all $\bfz\in\bfV_\sigma$.
Moreover, there is a non-negative, bounded and non-decreasing function  
$R_\mu(\cdot)$ on $[0,\infty)\times[0,\infty)$, depending only on 
$\mu >0$, such that
\begin{equation}\label{eq:5.6}
 |{\bfv}|^2_{W^{1,2}(0,T; {\bfH}_\sigma)}
   +\frac {\nu}2 \sup_{t \in [0,T]} |{\bfv}(t)|^2_{{\bfV}_\sigma} 
   \le R_\mu(|{\bfv}_0|_{{\bfV}_\sigma}, |\bfg|_{L^2(0,T; \bfH_\sigma})). 
\end{equation}
}

\noindent
For the solvability of $(H;u,{\bfv}_0,\bfg)^{\varepsilon}_{\mu}$ we apply
the general theory from Appendix III. 
To this end, we recall the following lemma, which is derived from 
the assumption
$u \in W^{1,2}(0,T;V^*_0)$ (hence $(\gamma_\mu u)^\varepsilon \in 
W^{1,2}(0,T; C(\overline \Omega)$ by \eqref{eq:2.8}). \vspace{0.5cm}

\noindent
{\bf Lemma 5.1 (cf. [11, Lemma 4.3] or [12, Lemma 2.2]).} 
{\it There exists a positive constant
$C_\mu$, depending only on $\mu$, which satisfies the following property:
for each $s, t \in [0,T]$ and
${\bfz}\in K_\mu((\gamma_\mu u)^\varepsilon;s)$ there is 
$\tilde {\bfz} \in K_\mu((\gamma_\mu u)^\varepsilon;t)$ such that
 $$ |\tilde {\bfz}-{\bfz}|_{{\bfH}_\sigma}\le C_\mu|(\gamma_\mu u)^\varepsilon
   (t)-(\gamma_\mu u)^\varepsilon(s)|_{C(\overline \Omega)},~~~
        \psi^t_\mu((\gamma_\mu u)^\varepsilon;\tilde{\bfz}) \le 
        \psi^s_\mu((\gamma_\mu u)^\varepsilon;{\bfz}). $$
}

Lemma 5.1 shows that problem \eqref{eq:5.4} 
can be handled in the general framework of Appendix with the set-up:
  $$X:= {\bfH}_\sigma, ~~\{\varphi^t(\cdot)\}:=
     \{\psi^t_\mu((\gamma_\mu u)^\varepsilon;\cdot)\} \in \Phi_c(M)~{\rm with~}
     M\geq |a|^2_{W^{1,2}(0,T)}$$
where
  $$ a(t):=C_\mu \int_0^t\left|\frac d{d\tau}(\gamma_\mu u)^{\varepsilon}(\tau)
   \right |_{C(\overline \Omega)} d\tau,~~b(\cdot)\equiv 0,  
    ~~h:={\bf{\mathcal G}}.$$
\vspace{0.5cm}

\noindent
{\bf Proof of Proposition 5.1.}
We observe that, cf.~Remark~5.1,  
   $$|{\bf{\mathcal G}}({\bfv},{\bfw})|_{{\bfH}_\sigma} \le
        p_0(\mu)|{\bfw}|_{{\bfV}_\sigma}, ~~\forall {\bfv} 
 \in K_{\mu}((\gamma_\mu u)^\varepsilon;t), \forall {\bfw} \in {\bfV}_\sigma 
$$ 
and
  $$ |({\bf{\mathcal G}}({\bfz}_1,{\bfz}_1)-{\bf{\mathcal G}}
({\bfz}_2,{\bfz}_2),
    {\bfz}_1-{\bfz}_2)_\sigma| \le 
     9 p_0 (\mu)|{\bfz}_1-{\bfz}_2|_{{\bfH}_\sigma}
  |{\bfz}_1-{\bfz}_2|_{{\bfV}_\sigma},$$
for all ${\bfz}_i\in K_{\mu}((\gamma_\mu u)^\varepsilon;t),~i=1,2$.
This shows
that the perturbation operator
  $$ \bfh(t,{\bfz}):={\bf{\mathcal G}}({\bfz},{\bfz}), ~\forall {\bfz} \in 
    K_\mu((\gamma_\mu u)^\varepsilon; t),~\forall t \in [0,T],$$
fulfills condition $(h4)$ in Appendix III. Also, it is easy to see that
this operator fulfills the other conditions 
$(h1)-(h3)$. 
Therefore, on account of Proposition III(1) in Appendix, 
the problem \eqref{eq:5.4}, namely $(H; u,{\bfv}_0,\bfg)^\varepsilon_\mu$, 
has one and only one solution ${\bfv} \in W^{1,2}(0,T;{\bfH}_\sigma)$
such that $t \to \psi^t_\mu((\gamma_\mu u)^\varepsilon;{\bfv}(t))$ is 
absolutely continuous on $[0,T]$.
This implies that ${\bfv} \in C([0,T];{\bfV}_\sigma)$. 
By Proposition III(2),
we obtain an estimate of the form \eqref{eq:5.6}.

Finally, we prove \eqref{eq:5.5}. Multiply
the inclusion in
 \eqref{eq:5.4}
by ${\bfv}$ and integrate in time over $[0,t]$  to get
  $$ \frac 12 |{\bfv}(t)|^2_{{\bfH}_\sigma} +\nu \int_0^t
      |{\bfv}(\tau)|^2_{{\bfV}_\sigma} d\tau +
 \int_0^t \langle {\bf{\mathcal G}}({\bfv}(\tau), {\bfv}(\tau)), {\bfv}(\tau)\rangle
  _\sigma d\tau
     \le \frac 12|{\bfv}_0|^2_{{\bfH}_\sigma}+\int_0^t(\bfg(\tau), \bfv(\tau))_\sigma d\tau.$$
By Remark~5.1(a), we immediately obtain \eqref{eq:5.5} 
from the above inequality. 
\hfill $\Box$
\vspace{0.5cm}

\noindent
{\bf Proposition 5.2.} {\it Let $\mu$ be any small positive number, and 
let $u_0 \in H$ with $0\le u_0 \le u^*$ a.e.~on 
$\Omega$.
Let $\{u_n\}$ be a bounded sequence in
$W^{1,2}(0,T;V^*_0)$ with $0\le u_n \le u^*$ a.e.~on $Q$ such that 
$u_n(0)=u_0$ and $ u_n \to u~~{\it in~}C([0,T];V^*_0)$.
Then, for any $\bfv_0\in K_\mu(\gamma_\mu u)^\eps;0)$ 
and any $\bfg\in L^2(0,T; \bfH_\sigma)$,
the solution ${\bfv}_n$ of $(H;u_n,{\bfv}_0,\bfg)^\varepsilon_\mu$ 
converges to the solution ${\bfv}$ 
of $(H;u,{\bfv}_0,\bfg)^\varepsilon_\mu$ in the sense that
  $${\bfv}_n \to {\bfv}~{\it in~}C([0,T];{\bfH}_\sigma)\cap 
    L^2(0,T; {\bfV}_\sigma)~
    {\it and~weakly ~in~}W^{1,2}(0,T;{\bfH}_\sigma).$$
}

\noindent
{\bf Proof.} Let us recall (cf.~\eqref{eq:2.8}) that 
$(\gamma_\mu u_n)^\varepsilon \to (\gamma_\mu u)^\varepsilon$ 
uniformly on $\overline Q$
(as $n\to \infty$). We show first that, for every $t \in [0,T]$,  
$\psi_\mu^t((\gamma_\mu u_n)^\varepsilon;\cdot) \to
\psi_\mu^t((\gamma_\mu u)^\varepsilon; \cdot)$ on ${\bfH}_\sigma$ 
in the sense of Mosco, as $n\to \infty$ (cf.~Appendix II). 
To this end, assume that 
$\{{\bfz}_n\}$ is any sequence
in ${\bfH}_\sigma$
with $\liminf_{n\to \infty}\psi_\mu^t((\gamma_\mu u_n)^\varepsilon;{\bfz}_n) 
< \infty$ and
${\bfz}_n \to {\bfz}$ weakly in ${\bfH}_\sigma$. It is enough to consider the
case ${\bfz}_n \in K_\mu((\gamma_\mu u_n)^\eps; t)$ and  
$\{{\bfz}_n\}$ is bounded in ${\bfV}_\sigma$.
In this case, 
$|{\bfz}_n(x)| \le p_\mu((\gamma_\mu u_n)^\varepsilon(x,t))$ 
for a.e.~$x\in \Omega$ and
${\bfz}_n \to {\bfz}$ in ${\bfH}_\sigma$ by the boundedness of 
$\{{\bfz}_n\}$ in ${\bfV}_\sigma$.This strong convergence yields 
$|{\bfz}(x)| \le p_\mu((\gamma_\mu u)^\varepsilon(x,t))$ 
for a.e.~$x \in \Omega$, 
so that
${\bfz} \in K_\mu((\gamma_\mu u)^\varepsilon;t)$, namely 
$\psi_\mu^t((\gamma_\mu u)^\varepsilon;{\bfz})<\infty$. As a consequence we 
have, as
${\bfz}_n \to {\bfz}$ weakly in ${\bfV}_\sigma$, that
  $$ \liminf_{n\to \infty}\psi_\mu^t((\gamma_\mu u_n)^\varepsilon;{\bfz}_n)
 \geq
   \psi_\mu^t((\gamma_\mu u)^\varepsilon;{\bfz}).$$
Next, let ${\bfz}$ be any function in $K_\mu((\gamma_\mu u)^\varepsilon;t)$. 
Then we construct the function ${\bfz}_n$ by: 
\begin{equation}\label{eq:5.7}
{\bfz}_n(x)=\left(1- \frac 1\mu|p_\mu(\gamma_\mu u_n)^\varepsilon(t))
    -p_\mu((\gamma_\mu u)^\varepsilon(t))|
    _{C(\overline \Omega)}\right){\bfz}(x),~~x \in \Omega. 
    \end{equation}
Since $(\gamma_\mu u_n)^\varepsilon \to (\gamma_\mu u)^\varepsilon$ uniformly 
on $\overline Q$ as $n \to \infty$ and 
$\frac {p_\mu((\gamma_\mu u)^\varepsilon)}{\mu} \geq 1$ by \eqref{eq:5.3},
it follows (cf. [12, Lemma~2.2]) that 
${\bfz}_n \in K_\mu((\gamma_\mu u_n)^\varepsilon;t)$ for all large $n$ and
${\bfz}_n \to {\bfz}$ in ${\bfV}_\sigma$ 
(hence $\psi_\mu^t((\gamma_\mu u_n)^\varepsilon;{\bfz}_n) \to 
  \psi_\mu^t((\gamma_\mu u)^\varepsilon;{\bfz}_n$)). 
  Accordingly,
$\psi_\mu^t((\gamma_\mu u_n)^\varepsilon;\cdot) \to
\psi_\mu^t((\gamma_\mu u)^\varepsilon; \cdot)$ on ${\bfH}_\sigma$ in the sense 
of Mosco.

We are now in a position to apply Proposition II 
to the sequence of problems
\begin{equation}\label{eq:5.8}
{\bfv}'_n(t) +\partial \psi_\mu^t((\gamma_\mu u_n)^\varepsilon; {\bfv}_n(t))
  +{\bf{\mathcal G}}({\bfv}_n(t),{\bfv}_n(t)) \ni \bfg(t)~{\rm in ~}{\bfH}_\sigma,~~
  {\bfv}_n(0)={\bfv}_0. 
\end{equation}
We note that all the families 
$\{\psi^t_\mu((\gamma_\mu u_n)^\varepsilon; \cdot)\}$, $n=1,2,\cdots$,
belong to the same class $\Phi_c(M)$ for a large number $M$,
since, by assumption and \eqref{eq:2.8}, 
$\{(\gamma_\mu u_n)^\eps\}$ is uniformly bounded in 
$W^{1,2}(0,T; C(\overline \Omega))$. 
Therefore, by virtue of Proposition 5.1, problem \eqref{eq:5.8} has
one and only one solution ${\bfv}_n$, and
the uniform estimates \eqref{eq:5.5} and \eqref{eq:5.6}
hold for each ${\bfv}_n$.
Hence, there is a subsequence $\{{\bfv}_{n_k}\}$ of $\{{\bfv}_n\}$ 
such that
\begin{equation}\label{eq:5.9}
{\bfv}_{n_k} \to {\bfv}~{\rm weakly~in~}
   W^{1,2}(0,T;{{\bfH}_\sigma})~
    {\rm and~weakly}^*~{\rm in~}L^\infty(0,T;{\bfV}_\sigma)
\end{equation}
(as $k \to \infty$), which implies that 
\begin{equation}\label{eq:5.10}
{\bfv}_{n_k} \to {\bfv}~\text{ in }
C([0,T];{\bfH}_\sigma)\  \text{ and } \ 
{\bf{\mathcal G}}({\bfv}_{n_k}, {\bfv}_{n_k}) \to {\bf{\mathcal G}}({\bfv},{\bfv}) \text{ weakly in }L^2(0,T;{\bfH}_\sigma).
\end{equation}
Therefore, by Proposition II, 
${\bfv}$ solves \eqref{eq:5.4}. 
Furthermore, by uniqueness of solution to \eqref{eq:5.4},
we obtain \eqref{eq:5.9} without extracting any subsequence from 
$\{{\bfv}_n\}$.

It remains to show the convergence ${\bfv}_n \to {\bfv}$ in 
$L^2(0,T;{\bfV}_\sigma)$. 
We consider the function $\tilde {\bfv}_n$ given by 
 $$\tilde {\bfv}_n(x,t):=\left(1-\frac 1\mu 
|p_\mu((\gamma_\mu u_n)^\varepsilon(t))-p_\mu((\gamma_\mu u)^\varepsilon(t))|_
{C(\overline \Omega)}
     \right){\bfv}(x,t),~~(x,t) \in Q.$$
Just as for \eqref{eq:5.7} above, we observe from [12, Lemma~2.2] again that 
\begin{equation}\label{eq:5.11}
\tilde {\bfv}_n(t) \in K_\mu((\gamma_\mu u_n)^\varepsilon;t)
\text{ for all large }n
\quad \text{ and }\quad 
\tilde{\bfv}_n \to {\bfv} \text{ in } L^2(0,T;{\bfV}_\sigma).
\end{equation}
Since $g-{\bfv}'_n-{\bf{\mathcal G}}({\bfv}_n,{\bfv}_n) \in 
\partial \psi_\mu^t((\gamma_\mu u_n)^\varepsilon;{\bfv}_n)$, 
it follows from \eqref{eq:5.sub} that
 \begin{equation}\label{eq:5.12} 
 \begin{array}{l}
\displaystyle {\int_0^T
\bigg({\bfv}'_n(t)+{\bf{\mathcal G}}\big({\bfv}_n(t),{\bfv}_n(t)\big)-\bfg(t),
      {\bfv}_n(t)-\tilde{\bfv}_n(t)\bigg)_\sigma dt}~~~~~~~~~~~~~~~~~~~~~\\[0.15cm]
 ~~~~~~~~~~~~~~~~~\displaystyle{ \le \nu \int_0^T \int_\Omega \nabla 
     {\bfv}_n(x,t)\cdot
     \nabla (\tilde {\bfv}_n(x,t)-{\bfv}_n(x,t))dxdt.}
  \end{array} 
  \end{equation}
Here, the left hand side of \eqref{eq:5.12} tends to $0$ 
as $n\to \infty$, since, from~\eqref{eq:5.10} and \eqref{eq:5.11},
 \begin{multline*}
 \liminf_{n \to \infty}
     \int_0^T
     \bigg({\bfv}'_n(t)+{\bf{\mathcal G}}\big({\bfv}_n(t),{\bfv}_n(t)\big)-\bfg(t),
      {\bfv}_n(t)-\tilde{\bfv}_n(t)\bigg)_\sigma dt 
 \\  \geq 
   \frac 12|{\bfv}(T)|^2_{{\bfH}_\sigma} 
    -\frac 12|{\bfv}_0|^2_{{\bfH}_\sigma}
     -\int_0^T({\bfv}'(t),{\bfv}(t))_\sigma dt =0.
     \end{multline*}
Therefore
  $$ \limsup_{n \to \infty}\int_0^T \int_\Omega|\nabla {\bfv}_n|^2dxdt
      \le \lim_{n \to \infty}\int_0^T\int_\Omega \nabla {\bfv}_n\cdot 
        \nabla \tilde{\bfv}_n dxdt
      =\int_0^T\int_\Omega|\nabla {\bfv}|^2dxdt.$$
This implies  ${\bfv}_n \to {\bfv}$ in $L^2(0,T;{\bfV}_\sigma)$. 
\hfill $\Box$

\section{Approximate full system and its convergence}\label{sec:full}
Let $\varepsilon$ be a small positive parameter and fix it. For each small
$\mu>0$, consider the coupling
$P^{\varepsilon}_\mu :=\{(B;w,{\bfv},u_0)_\mu,
(N;u,{\bfv},w_0)_\mu, (H;u,{\bfv}_0,\bfg)^\varepsilon_\mu\}$ as the approximation
to our problem $P^{\varepsilon}=\{(B;w,{\bfv},u_0), (N;u,{\bfv},w_0),
(H;u,{\bfv}_0,\bfg)^{\varepsilon}\}$. \vspace{0.5cm}

More precisely, a triplet $\{u_\mu,w_\mu, {\bfv}_\mu\}$ is called a solution of
$P^{\varepsilon}_\mu$, if 
\begin{description}
\item{(a)} $u_\mu \in W^{1,2}(0,T;V^*_0)$,
{$t \to |\hat \beta(u_\mu(t))|_{L^1(\Omega)}$
is absolutely continuous on $[0,T]$},
and $u_\mu$
is the solution of 
$(B;w_\mu,{\bfv}_\mu, u_0)_\mu$;
\item{(b)} $w_\mu \in L^2(0,T;V),~w'_\mu \in
L^2(0,T;V^*)$, $0 \le w_\mu \le 1$ a.e.~on 
$Q$ and $w_\mu$ is the solution of 
$(N;u_\mu,{\bfv}_\mu, w_0)_\mu=(N; \rho_\mu*u_\mu, \bfv_\mu, w_0)$;
\item{(c)} ${\bfv}_\mu \in W^{1,2}(0,T;{\bfH}_\sigma) \cap C([0,T];
{\bfV}_\sigma)$ and ${\bfv}_\mu$ is the solution of 
$(H;u_\mu,{\bfv}_0,\bfg)^\varepsilon_\mu$.
\end{description}
\vspace{0.25cm}

\noindent
{\bf Theorem 6.1.} {\it Let $\mu \in(0,\delta_0)\cap(0,1)$ 
with $\mu< p_0(\mu)$. Assume that $u_0 \in H$ is such that 
$\hat \beta(u_0)\in L^1(\Omega), ~w_0\in H$
with $0\le w_0\le 1$ a.e.~on $\Omega$ and
${\bfv}_0 \in {\bfV}_\sigma \cap C(\overline\Omega)^3$ with
$|{\bfv}_0| < p_0(u_0^\varepsilon)$ on $\overline \Omega$,
where $u^\varepsilon_0(x)= \int_\Omega \rho_\varepsilon(x-y)u_0(y)dy$ for
all $x \in \overline \Omega$. Let $\bfg\in L^2(0,T;\bfH_\sigma)$.
Then, for all small positive number $\mu$
the approximate system $P^{\varepsilon}_\mu$
has at least one solution $\{u_\mu, u_\mu, {\bfv}_\mu\}$.}\vspace{0.2cm}

\noindent
{\bf Proof.} 
We put
   $$ X(u_0):=\left\{u~\left|~
        \begin{array}{l}
       \displaystyle{  | u|^2_{W^{1,2}(0,T;V^*_0)} + 
      \sup_{t \in [0,T]}|\hat\beta(u(t))|_{L^1(\Omega)}} \\[0.3cm]
    \displaystyle{ \le B_0\left(T^{\frac 12} 
|{\bfv}_0|_{{\bfH}_\sigma},
        |\hat\beta(u_0)|_{L^1(\Omega)}\right),}\\[0.3cm]
        u(0)=u_0
      \end{array} \right. \right \}, $$
where $B_0(\cdot)$ is the same function 
 as in \eqref{eq:3.6} of Proposition 3.1. Note that
$X(u_0)$ is non-empty, compact and convex in $C([0,T];V^*_0)$.
By assumption, for each $u \in X(u_0)$ we see that 
$|{\bfv}_0| \le p_\mu((\gamma_\mu u)^\varepsilon(\cdot,0))$ on $\Omega$ 
for all small $\mu>0$, 
since $(\gamma_\mu u)^\varepsilon \to u^\varepsilon$
in $C(\overline Q)$ by $ \gamma_\mu u_0 \to u_0$ in $H$ as
$\mu \downarrow 0$. This implies that ${\bfv}_0 \in K_\mu((\gamma_\mu u)
^\varepsilon; 0)$ for all small $\mu >0$, so that 
$(H;u,{\bfv}_0,\bfg)^{\varepsilon}_\mu$ is uniquely solved. Now, denote
the solution by ${\cal S}_1u=:{\bfv}$.  Then, according to Proposition 5.1,
${\bfv} \in W^{1,2}(0,T,{\bfH}_\sigma)\cap 
C([0,T];{\bfV}_\sigma)$ and there is a positive constant 
$R_\mu(|{\bfv}_0|_{{\bfV}_\sigma},|\bfg|_{L^2(0,T;\bfH_\sigma)})=:R_\mu$, 
depending on the parameter $\mu$, $|{\bfv}_0|_{{\bfV}_\sigma}$ 
and $|\bfg|_{L^2(0,T;\bfH_\sigma)}$, such that (cf. \eqref{eq:5.6})
   $$ |{\bfv}|^2_{W^{1,2}(0,T;{\bfH}_\sigma)} + 
    \frac \nu 2 \sup_{t \in[0,T]}|{\bfv}(t)|^2_{{\bfV}_\sigma} \le R_\mu.$$
Put 
  $$Y({\bfv}_0):=\left\{{\bfv}~\left|~|{\bfv}|^2_{W^{1,2}(0,T;{\bfH}_\sigma)} +
    \frac \nu 2 \sup_{t \in[0,T]}|{\bfv}(t)|^2_{{\bfV}_\sigma} \le R_\mu,
     ~{\bfv}(0)={\bfv}_0 \right.\right\}.$$
Then, ${\cal S}_1$ is a mapping from $X(u_0)$ into $Y({\bfv}_0)$.

Next, for each pair of $u \in X(u_0)$ and ${\bfv}$ of $Y({\bfv}_0)$ we solve
$(N;u,{\bfv},w_0)_\mu$ and denote its 
solution by ${\cal S}_2(u,{\bfv})=: w$. By Proposition 4.1, 
estimate \eqref{eq:4.5} holds for $w$. Put
  $$ Z(w_0):=\{w\in L^2(0,T;V)\cap L^\infty
 (0,T;H)~|~\eqref{eq:4.5} ~{\rm holds},~w(0)=w_0 \}.$$
Then ${\cal S}_2$ is a mapping from 
$X(u_0)\times Y({\bfv}_0)$ into $Z(w_0)$.

Furthermore, for each $w \in Z(w_0)$ and ${\bfv} \in Y({\bfv}_0)$ 
we solve $(B;w,{\bfv},u_0)_\mu$
and denote the solution by ${\cal S}_3(w,{\bfv})=: \bar u$. It 
follows from Proposition 3.1 with \eqref{eq:3.6}
that $\bar u \in X(u_0)$,
so that ${\cal S}_3$ can be considered as a mapping from 
$Z(w_0)\times Y({\bfv}_0)$ into 
$X(u_0)$. 
Finally we define a mapping ${\cal S}$ from $X(u_0)$ into itself
by
    $$ {\cal S}u:= {\cal S}_3\left({\cal S}_2(u,{\cal S}_1u), {\cal S}_1(u))
    \right),
      ~~\forall u \in X(u_0).$$

In order to apply the fixed-point theorem for
compact mappings we show continuity of 
${\cal S}$.
Assume that $u_n \in X(u_0)$ and $u_n \to u$ in $C([0,T];V^*_0)$. Then,
by the definition of $X(u_0)$, 
$u_n \to u$ weakly in $W^{1,2}(0,T;V^*_0)$. 
Moreover, by lower semicontinuity of $\hat{\beta}$, 
$u \in X(u_0)$.
Therefore, $(\gamma_\mu u_n)^\varepsilon:=\rho_\varepsilon*(\gamma_\mu u_n)
\in W^{1,2}(0,T;C(\overline \Omega))$ by \eqref{eq:2.8} and 
$$|(\gamma_\mu u_n)^\varepsilon|_{W^{1,2}(0,T;C(\overline \Omega))} \le 
R'_\mu,~~
   |(\gamma_\mu u)^\varepsilon|_{W^{1,2}(0,T;C(\overline \Omega))} \le R'_\mu,$$
for some positive constant $R'_\mu$ and 
$ (\gamma_\mu u_n)^\varepsilon \to (\gamma_\mu u)^\varepsilon
 ~{\rm uniformly ~on~}\bar{Q}$ (as $n \to \infty$).
As we have seen in Propositions 5.1 and 5.2, the problem 
$(H;u_n,{\bfv}_0,\bfg)^{\varepsilon}_\mu$ 
has a unique solution 
${\bfv}_n$ in $W^{1,2}(0,T;{\bfH}_\sigma)\cap C([0,T];{\bfV}_\sigma)$. Moreover, 
${\bfv}_n$
converges in $C([0,T];{\bfH}_\sigma)$
$\cap L^2(0,T;{\bfV}_\sigma)$
and weakly in 
$W^{1,2}(0,T;{\bfH}_\sigma)$ to
the solution ${\bfv}$ of $(H;u,{\bfv}_0,\bfg)^{\varepsilon}_\mu$. 
This fact implies that 
 $${\cal S}_1u_n \to {\cal S}_1u~{\rm in~} 
C([0,T];{\bfH}_\sigma)\cap L^2(0,T; {\bfV}_\sigma)~{\rm and
~weakly~ in~} W^{1,2}(0,T;{\bfH}_\sigma).$$

Next, as for the sequence $\{{\cal S}_2(u_n,{\bfv}_n)\}$ with 
${\bfv}_n:={\cal S}_1u_n$, we obtain from Proposition~4.1 that
$(N;u_n,{\bfv}_n,w_0)_\mu$
has a unique solution 
$w_n:={\cal S}_2(u_n,{\bfv}_n)$
in 
$L^2(0,T;V)$ with $w'_n \in L^2(0,T;V^*)$ 
and $0\le w_n\le 1$ a.e.~$Q$, 
satisfying the uniform estimate
\begin{equation}\label{eq:6.1}
|w_n(t)|^2_H+2c_d\int_0^t |\nabla w_n|^2_H
  d\tau \le e^{2u^*L(f)T} |w_0|^2_H,~~
  \forall t \in [0,T].  
\end{equation}
This estimate shows that $w'_n=-\partial_{V^*}\Phi^t(\rho_{\mu}*u_n; w_n)-{\bfv}_n\cdot
\nabla w_n-f(w_n)\rho_{\mu}*u_n$ is bounded
in $L^2(0,T;V^*)$ (as $n\to \infty$), so that
it follows from the Aubin's compactness theorem [2]
that $\{w_n\}$ is relatively compact in 
$L^2(Q)$. Now, we choose a subsequence 
$\{w_{n_k}\}$ of $\{w_n\}$ so as to satisfy $w_{n_k}\to w$ 
in $L^2(Q)$ (as $k\to \infty$) for some
function $w$. Then, by \eqref{eq:6.1}, 
$w_{n_k} \to w$ weakly in $L^2(0,T;V)$,
$w'_{n_k}\to w'$
weakly in $L^2(0,T;V^*)$ and $0\le w\le 1$
a.e.~on $Q$.  
Besides, since $\rho_\mu*u_{n_k} \to \rho_\mu*u$ in $L^2(Q)$ 
(cf. \eqref{eq:2.9}), we have
 $$f^*_{n_k}:={\bfv}_{n_k}\cdot \nabla w_{n_k}+f(w_{n_k})\rho_\mu*u_{n_k}
\to f^*:={\bfv}\cdot \nabla w +f(w)\rho_\mu*u~~
    {\rm in~} L^2(0,T;V^*).$$ 
As a consequence, letting $k \to \infty$,
we see that $w$ is the solution of 
$(N; u,{\bfv},w_0)_\mu$. 
Since the solution of $(N; u,{\bfv},w_0)_\mu$ 
is unique, the above convergences hold without
extracting any subsequence from $\{w_n\}$, that is, 
$${\cal S}_2(u_n,{\bfv}_n) = w_n \to w:={\cal S}_2(u,{\bfv})~
{\rm in~} L^2(Q) \text{ and weakly in }L^2(0,T;V),
$$
$$
 \text{ and }w'_n \to w'\text{ weakly in } L^2(0,T;V^*).
$$ 
Moreover, by the convergences of $\{{\bfv}_n\}$ and $\{w_n\}$ obtained above,
Proposition 3.2 implies that
the solution of $(B;w_n,{\bfv}_n,u_0)_\mu $ converges to that of 
$(B;w,{\bfv}, u_0)_\mu $ in $C([0,T];V^*_0)$.
Namely, 
$${\cal S}u_n={\cal S}_3(w_n,{\bfv}_n) \to {\cal S}_3(w,{\bfv})
={\cal S}u \text{ in } C([0,T];V^*_0).
$$
Thus, ${\cal S}$ is continuous in 
$C([0,T];V^*_0)$. 
Accordingly, it follows from the Schauder fixed-point theorem
that ${\cal S}$ admits at least one fixed-point, $u_\mu={\cal S}u_\mu$. 
It is easy to see that this fixed-point $u_\mu$ with the solutions 
${\bfv}_\mu$ of 
$(H;u_\mu,{\bfv}_0,\bfg)^\varepsilon_\mu$ and 
$w_\mu$ of $(N;u_\mu,{\bfv}_\mu,w_0)_\mu$ gives
 a set of solutions of
our problem $P^{\varepsilon}_\mu$. \hfill $\Box$
\vspace{0.5cm}

Now, we summarize the uniform estimate on approximate solutions
$\{u_\mu,w_\mu,{\bfv}_\mu\}$; we have automatically
\begin{equation}\label{eq:6.2}
 0\le u_\mu \le u^*,~~0\le w_\mu \le 1~~
    \text{ a.e.~in~}Q. 
\end{equation}
Furthermore, by our construction of approximate solutions, 
there is a positive constant $A_0$ depending only 
on the data $u_0,~w_0,~{\bfv}_0,~\bfg,~\beta,~f,~d$ and $p_0$ such that
\begin{multline}\label{eq:6.3}
   |u_\mu|_{W^{1,2}(0,T;V^*_0)} 
        + \max_{t \in [0,T]}|\hat \beta(u_\mu(t))|_{L^1(\Omega)}
   + |w_\mu|_{L^2(0,T;V)}\\[0.12cm]
   +|w'_\mu|_{L^2(0,T;V^*)}+|{\bfv}_\mu|_{L^2(0,T;{\bfV}_\sigma)}
   +|{\bfv}_\mu|_{L^\infty(0,T;{\bfH}_\sigma)} \le A_0 
\end{multline}
for all small $\mu>0$. 
%
On account of the uniform estimates \eqref{eq:6.2}, \eqref{eq:6.3}, 
Lemma 3.1 and Proposition 4.2, there is a sequence
$\{\mu_n\}$ with $\mu_n \downarrow 0$ (as $n\to \infty$) and a triplet
$\{u,w,{\bfv}\}$ of functions such that
\begin{equation}\label{eq:6.4}
u_n:=u_{\mu_n} \to u~{\rm in~}L^2(Q)~{\rm and ~weakly~in~}W^{1,2}(0,T;V^*_0)
\end{equation}
\begin{equation}\label{eq:6.5}
w_n:=w_{\mu_n} \to w~{\rm in~}L^2(Q),
~{\rm weakly~in~}L^2(0,T;V),~w'_n \to w'~{\rm weakly~in~}L^2(0,T;V^*),
\end{equation}
\begin{equation}\label{eq:6.6}
{\bfv}_n:= {\bfv}_{\mu_n} \to {\bfv}~{\rm weakly~in~}L^2(0,T;{\bfV}_\sigma)~
    \text{ and~weakly}^*~{\rm in~}L^\infty(0,T;{\bfH}_\sigma).
\end{equation}

In the rest of this section we shall show that $\{u,w,{\bfv}\}$ is a
solution of the limit problem $P^\varepsilon$. 
To this end, we make use of recent important results about the convergence 
of $\{\bfv_n\}$, which was obtained in the authors' work [12].
\vspace{0.5cm}

\noindent
{\bf Theorem 6.2.} {\it Let $\varepsilon$ be a small positive number and fix it. Assume that $u_0 \in H$ with
$\hat \beta(u_0) \in L^1(\Omega)$, $w_0 \in H$ with
$0 \le w_0 \le 1$ a.e.~on $\Omega$, $\bfg\in L^2(0,T;\bfH_\sigma)$, and 
\begin{equation}\label{eq:6.7}
{\bfv}_0 \in {\bfW}_{0,\sigma}^{1,4}(\Omega),
 ~{\rm supp}({\bfv}_0) \subset \{x \in \Omega~|~
 p_0(u^\varepsilon_0(x))>0\},~|{\bfv}_0| <
 p_0(u^\varepsilon_0)~{\it on~}\overline \Omega,
\end{equation}
where $u^\varepsilon_0=\rho_\varepsilon*u_0$. 
Then there exists at least one set of functions
$\{u, w, {\bfv}\}$ such that
\begin{description}
\item{(i)} $u$
is a solution of $(B;w,{\bfv},u_0)$ in the 
sense of Definition 3.1.
\item{(ii)} $w$ is a solution of 
$(N;u,{\bfv},w_0)$ in the sense of Definition 4.1.
\item{(iii)} ${\bfv}$ is a weak solution of
$(H;u,{\bfv}_0,\bfg)^\varepsilon$ in the sense of Definition 5.1. 
\end{description}
}

\noindent
{\bf Proof.} Fix $\varepsilon \in (0,1]$.
Let $\{u_n, w_n, {\bfv}_n\}$ be the same sequence of approximate solutions
as in \eqref{eq:6.4}--\eqref{eq:6.7} with the limit $\{u,w,{\bfv}\}$. 
As for the convergences of $(B;w_n,{\bfv}_n,u_0)_{\mu_n}$ and 
$(N; u_n,{\bfv}_n,w_0)_{\mu_n}$,
by \eqref{eq:6.4} and \eqref{eq:6.5}, we see that
$$ \rho_{\mu_n}*(\gamma_{\mu_n}u_n) \to u,
\quad 
\rho_{\mu_n}*u_n \to u
\quad\text{ in }L^2(Q),
$$
and 
$$
\rho_{\mu_n}*w_n \to w~\text{ in~}L^2(Q)\text{ and
weakly in }L^2(0,T;V).
$$
Therefore, by Propositions 3.2 and 4.2 the limits $u$ and $w$ 
are solutions  of $(B;w,{\bfv},u_0)$ and $(N; u,{\bfv},w_0)$, 
respectively. Thus {\it (i)} and {\it (ii)} hold. 
In order to complete the proof of Theorem 6.2 it remains to
prove {\it (iii)}.
 
Actually, we are going to show that {\it (iii)} 
is a direct consequence of [12], putting
   $$ p(x,t):=p_0(u^\eps(x,t)),~~p_n(x,t):=p_{\mu_n}((\gamma_{\mu_n}u_n)^\eps
(x,t)),
    ~~\forall (x,t) \in \overline Q,~~\forall n \in {\bf N}.$$
From this definition of $p,~p_n$ and the fact 
$(\gamma_{\mu_n}u_n)^\eps \to u^\eps$ in
$C(\overline Q)$ it is easy to see that
  \begin{equation}\label{eq:6.8}
  \left\{ 
      \begin{array}{l}
      0< p_n < \infty ~{\rm on ~}\overline Q,~\forall n \in {\bf N},\\
      \forall \kappa \in (0,\infty),~p_n \to p~{\rm uniformly~on~}
     \{(x,t) \in \overline Q~|~p(x,t)\le \kappa\}
      \end{array} \right. 
  \end{equation}
and
 \begin{equation}\label{eq:6.9}
  \left\{
   \begin{array}{l}
     \forall M>0~{\rm sufficiently~large~},
     ~\exists n_M \in {\bf N}~{\rm such~that~}\\
    ~~~~~
    p_n > M~~{\rm on~} \{(x,t)\in \overline Q~|~p(x,t) > M\},
    ~~\forall n\geq n_M.
    \end{array} \right. 
    \end{equation}
Under \eqref{eq:6.8} and \eqref{eq:6.9} it is proved in [12; Lemma 4.1, Theorem 1.1] 
that the  sequence of solutions $\bfv_n$ of variational inequalities 
$(H;u_n, \bfv_0, \bfg)^\eps_{\mu_n}$ 
of the Navier-Stokes type converges to $\bfv$, and,
in addition to \eqref{eq:6.6}:
\begin{description}
\item{(1)} For every $t \in [0,T]$, $\bfv_n(t) \to \bfv(t)$ weakly in 
$\bfH_\sigma$ and $|\bfv(x,t)|\le p_0(u^\varepsilon(x,t))$ 
for a.e.~$x \in \Omega$.
\item{(2)} $\bfv_n \to \bfv$ (strongly) in $L^2(0,T;{\bfH}_\sigma)$
\item{M}oreover: 
\item{(3)} For any test function 
${\boldsymbol{\mathcal \eta}} \in {\boldsymbol{\mathcal K}}(u^\eps)$, 
given in Definition~5.1, the real-valued function 
$t \to (\bfv(t),\boldsymbol{\mathcal \eta}(t))_\sigma$ 
is of bounded variation on $[0,T]$.
\item{(4)} The limit ${\bfv}$ satisfies \eqref{eq:5.2}.
\end{description}
In fact, by virtue of (1) and (2), the nonlinear term 
$(\bfv_n\cdot \nabla)\bfv_n$ converges to 
$(\bfv\cdot \nabla)\bfv$ in $L^{\frac 43}(0,T;W^{-1,\frac 43}(\Omega)^3)$ (the 
dual space of 
$L^4(0,T; W^{1,4}_0(\Omega)^3)$). Hence 
we can arrive at the variational
inequality \eqref{eq:5.2} by  integrating by parts in time and
letting $n\to \infty$ in the variational inequality equivalent to
$(H; u_{\mu_n},\bfv_0,\bfg)^\eps_{\mu_n}$. 
For the detailed proof, see [12]. 
\hfill $\Box$ \vspace{0.25cm}


\noindent
{\bf Remark 6.1.} In Theorem 6.2, we do not know whether 
${\bfv}(t)$ is continuous in time or not. However the initial condition
${\bfv}(0)={\bfv}_0$ makes sense, because ${\bfv}(t)$ is defined for
every $t \in [0,T]$ and the real-valued function 
$t \to ({\bfv}(t),{\boldsymbol{\mathcal \eta}}(t))_\sigma$ is of bounded 
variation on $[0,T]$ for any test function $\boldsymbol{\mathcal \eta}$.
In particular, if supp$(\boldsymbol{\mathcal \eta})$ is included in the liquid 
region (namely the 
interior of $\{(x,t) \in \overline Q~|~u^\varepsilon(x,t)=0\}$), 
$({\bfv}(t), \boldsymbol{\mathcal \eta}(t))_\sigma$ is absolutely coninuous in
 $t$ on $[0,T]$ and $({\bfv}(0), \boldsymbol{\mathcal \eta}(0))_\sigma=
({\bfv}_0, \boldsymbol{\mathcal \eta}(0))_\sigma$.
For this result, see [12; Corollary 3.2, Remark 4.1].
\vspace{0.25cm}

\noindent
{\bf Remark 6.2.}
A number of open questions concerning the mathematical modeling
of biomass development remain. For instance, the limit problem 
as $\delta_0 \to 0$
is the sharp interface model mentioned in the introduction.
It is expected that this question
will be affirmatively solved. 
Another question is to characterize the limit 
procedure of $\varepsilon \to 0$; when the convolution 
parameter $\varepsilon$ tends to 0, 
in which class of evolution
inclusions the limit problem can be handled. 
This seems a very difficult question.

\section*{Appendix}
Let $X$ be a real Hilbert space with inner product $(\cdot,\cdot)_X$ and norm 
$|\cdot|_X$. 
For a fixed (large) positive number $M$ we denote by $\Phi(M)$ the set of
all families $\{\varphi^t(\cdot)\}_{t \in [0,T]}$ of non-negative proper, 
l.s.c. and convex functions $\varphi^t(\cdot)$ on $X$ satisfying the following
conditions $(\Phi1)$ and $(\Phi2)$:
\begin{description}
\item{($\Phi1$)} $\displaystyle{\min_{z \in X}\{|z|^2_X+\varphi^t(z)\} \le M}$
for all $t \in [0,T]$.\vspace{0.3cm}
\item{($\Phi2$)} There are non-negative real-valued functions 
$a(\cdot) \in W^{1,2}(0,T)$ and
$b(\cdot) \in W^{1,1}(0,T)$ satisfying 
  $$|a|^2_{W^{1,2}(0,T)}+|b|_{W^{1,1}(0,T)} \le M $$
and the following property that for each 
$s,~t \in [0,T]$ and $z \in D(\varphi^s)$ there is an element 
$\tilde z \in D(\varphi^t)$ such that
 $$|\tilde z -z|_X \le |a(t)-a(s)|(1+\varphi^s(z)^{\frac 12}), $$
 $$\varphi^t(\tilde z) -\varphi^s(z) \le |b(t)-b(s)|(1+\varphi^s(z)). $$
\end{description}
\vspace{0.5cm}

We recall the fundamental results (cf. [4, 13, 18]) on
the Cauchy problem
$$ CP(\varphi^t; f,u_0) ~~\left\{
    \begin{array}{l}
     u'(t)+\partial_X \varphi^t(u(t)) \ni f(t)~~{\rm in~}X, ~t \in [0,T],
     \\[0.3cm]
     u(0)=u_0,
    \end{array} \right. $$
where $ f \in L^2(0,T; X)$ and $u_0 \in D(\varphi^0)$ are prescribed as data. 
It is said that $u: [0,T] \to X$ is a (strong) solution of
$CP(\varphi^t;f,u_0)$, if $u \in W^{1,2}(0,T;X)$, $u(0)=u_0$ and $f(t)-u'(t) \in \partial_X \varphi^t(u(t))$
for a.e.~$t \in (0,T)$, where $u'(t):=\frac {d u(t)}{dt}$.
\vspace{0.2cm}

We denote by $\Phi_c(M)$ the subclass of all families $\{\varphi^t\}$ in 
$\Phi(M)$ that satisfy the condition of level set compactness:
 $$ \{z \in X~|~\varphi^t(z)\le r\}~{\rm is~ cpmpact~ in~} X,~\forall r>0,~
    \forall t\in [0,T].$$

\noindent
{\bf [I] Existence and uniqueness}\vspace{0.15cm}

First of all, we recall the results on the existence, uniqueness and uniform 
estimates of solutions upon data for 
$CP(\varphi^t;f,u_0)$.
\vspace{0.3cm}

\noindent
{\bf Proposition I. (cf. [13; Chapter 1])}. 
{\it Let $\{\varphi^t\}\in \Phi(M)$. 
Then we have:
\begin{description}
\item{(1)} For each $f \in L^2(0,T;X)$ and $u_0 \in D(\varphi^0)$
the Cauchy problem $CP(\varphi^t;f,u_0)$ admits one and only one solution $u$ 
such that
$ u \in W^{1,2}(0,T;X)$
and the function $t \to \varphi^t(u(t))$ is absolutely continuous on $[0,T]$. 
\item{(2)} Let $\{f_i,u_{i0}\} \in L^2(0,T;X)\times D(\varphi^0),
~i=1,2,$ be two sets of data and denote by $u_i$ the solution of 
$CP(\varphi^t;f_i, u_{i0})$. Then we have, for all $s,~t \in [0,T]$ with $s\le t$:
$$
   \frac 12 |u_1(t)-u_2(t)|^2_X \le \frac 12 |u_1(s)-u_2(s)|^2_X
   +\int_s^t (f_1(\tau)-f_2(\tau), u_1(\tau)-u_2(\tau))_X d\tau. $$
\item{(3)} There is a non-negative and non-decreasing function $A_1:= 
A_1(M; n_1, n_2,n_3): {\bf R}^4_+ \to {\bf R}_+$ such that
 $$ |u|^2_{W^{1,2}(0,T;X)}+ \sup_{0\le t\le T}\varphi^t(u(t))
     \le A_1(M;n_1,n_2,n_3), $$
as long as $u$ is the solution of $CP(\varphi;f,u_0)$ with 
$|f|_{L^2(0,T;X)}\le n_1$, $|u_0|_X \le n_2$ and $\varphi^0(u_0)\le n_3$.
\end{description} }

\noindent
{\bf [II] Convergence results}\vspace{0.15cm}

Next, we recall the concept of Mosco convergence (cf. [16]).
Let $\{\varphi_n\}$ a sequence of non-negative proper, l.s.c. and convex
function on $X$. Then it is said that $\{\varphi_n\}$ converges 
to a non-negative, proper l.s.c.~and convex function $\varphi$ on $X$
(as $n\to \infty$) in the sense of Mosco, if
the following two conditions $(m_1)$ and $(m_2)$ are satisfied:
\begin{description}
\item{$(m_1)$} If $z_n \to z$ weakly in $X$, then 
   $\quad\liminf_{n \to \infty}\varphi_n(z_n) \geq \varphi(z).$
\item{$(m_2)$} For every $z \in D(\varphi)$ there is a sequence $\{z_n\}$
in $X$ such that 
$$ z_n \to z~~{\rm in~}X,~~\varphi_n(z_n) \to \varphi(z).$$
\end{description}
For other  characterizations of the Mosco convergence
see e.g.~[1; Chapter 3], [14; section~8].
\vspace{0.25cm}

\noindent
{\bf Proposition II. (cf. [13; Theorem 2.7.1])} {\it Let $\{\varphi^t_n\}$ 
be a sequence
of families in $\Phi(M)$ and $\{\varphi^t\} \in \Phi(M)$ such that
$\varphi^t_n$ converges to $\varphi^t$ in the sense of Mosco on $X$ for 
every $t \in [0,T]$. Also, let $\{f_n\}$ be a sequence in $L^2(0,T;X)$
such that $f_n \to f$ in $L^2(0,T; X)$, 
and $\{u_{n0}\}$ be a sequence in $X$ such that $u_{n0} \in D(\varphi^0_n)$,
$\sup_{n \in{\bf N}}{\varphi^0_n(u_{n0})} <\infty$
and $u_{n0} \to u_0$ in $X$.
Then the solution $u_n$ of $CP(\varphi^t_n;f_n,u_{n0})$ converges to the
solution $u$
of $CP(\varphi^t; f,u_0)$ in the sense that
 $$ u_n \to u~{\it in~}C([0,T];X),~~\int_0^T \varphi^t_n(u_n(t))dt
  \to \int_0^T \varphi^t(u(t))dt $$
and 
    $$ u_n \to u~~{\it weakly~in~}W^{1,2}(0, T;X).$$
In particular, if $\{\varphi^t_n\} \in \Phi_c(M)$ and 
$\{\varphi^t\} \in \Phi_c(M)$,
then the condition  ^^ ^^ $f_n \to f$ in $L^2(0,T;X)$" is replaced 
by ^^ ^^ $f_n \to f$ weakly in $L^2(0,T;X)$"}
\vspace{0.5cm}

\noindent
{\bf [III] A perturbation result}\vspace{0.15cm}

Finally, we present a perturbation result.
Let $\{\varphi^t\} \in \Phi_c(M)$ and  
let $h(t,\cdot)$ be a single-valued mapping 
from $D(\varphi^t)$ into $X$ for each $t \in [0,T]$ such that
\begin{description}
\item{(h1)} if $v \in L^2(0,T;X)$ with $v(t) \in D(\varphi^t)$
for a.e.~$t \in [0,T]$, then $h(\cdot,v(\cdot))$ is strongly measurable
on $[0,T]$,
\item{(h2)} there are positive constants $\alpha_1,~\alpha_2,~\alpha_3$
such that
 $$ |h(t,z)|_X^2 \le \alpha_1\varphi^t(z)+\alpha_1|z|_X^2+\alpha_3,
     ~\forall t \in [0,T],~\forall z \in D(\varphi^t),$$
\item{(h3)} (demi-closedness) if $t_n \in [0,T]$, $z_n \in X$,
$\{\varphi^{t_n}(z_n)\}$ is bounded, $z_n \to z$ in $X$ and $t_n \to t$
(as $n\to \infty$),
then $h(t_n,z_n) \to h(t,z)$ weakly in $X$.
\item{(h4)} for each $\delta >0$ there exists a positive constant 
$C_\delta$ such that
   $$|(h(t,z_1)-h(t,z_2), z_1-z_2)_X| \le \delta(z^*_1-z^*_2,z_1-z_2)_X
    +C_\delta |z_1-z_2|^2_X, $$
   $$ \forall t \in [0,T], ~\forall z_i \in D(\partial_X\varphi^t),~
      \forall z^*_i \in \partial_X\varphi^t(z_i),~i=1,2. $$
\end{description}
\vspace{0.2cm}

Now, given
$f\in L^2(0,T;X)$ and $u_0\in D(\varphi^0)$, we consider
the following perturbation problem, denoted by $CP(\varphi^t, h;f,u_0)$,
 $$ CP(\varphi^t, h;f,u_0)~~\left \{
     \begin{array}{l}
      u'(t)+\partial_X \varphi^t(u(t))+h(t,u(t)) \ni f(t)~~{\rm in~}X,~
      t \in [0,T],\\[0.3cm]
      u(0)=u_0.
     \end{array} \right. $$
It is said that $u$ is a solution of $CP(\varphi^t, h;f,u_0)$, if it is
a solution of $CP(\varphi^t; f-h(\cdot,u(\cdot)),u_0)$, namely
if $u \in W^{1,2}(0,T;X)$, $u(0)=u_0$ and $f(t)-u'(t)-h(t,u(t)) 
\in \partial_X \varphi^t(u(t))$
for a.e.~$t \in (0,T)$.
As to this perturbation problem we have similar results
to Proposition I. \vspace{0.5cm}

\noindent
{\bf Proposition III. [18; Theorem 2.1]} {\it Let 
$\{\varphi^t\} \in \Phi_c(M)$ and
$h(\cdot,\cdot)$ be a single-valued mapping satisfying $(h1) - (h4)$.
Then we have:
\begin{description}
\item{(1)}
For each $f \in L^2(0,T;X)$ and $u_0 \in D(\varphi^0)$
problem $CP(\varphi,h; f,u_0)$ admits one and only one solution $u$ such that
$ u \in W^{1,2}(0,T;X)$ and 
the function $t \to \varphi^t(u(t))$ is absolutely continuous on $[0,T]$. 
\item{(2)} There is a non-negative and non-decreasing function 
$A_2:=A_2(M,h; n_1, n_2, n_3) :{\bf R}^3_+ \to {\bf R}_+$,
depending only on the class $\Phi_c(M)$, $h$ and three given positive 
constants $n_1,~n_2,~n_3$,
such that
 $$ |u|^2_{W^{1,2}(0,T;X)}+ \sup_{0\le t\le T} \varphi^t(u(t))
     \le A_2(M,h;n_1,n_2,n_3), $$
as long as $u$ is the solution of $CP(\varphi^t,h;f,u_0)$ with 
$|f|_{L^2(0,T;X)}\le n_1$ and $u_0 \in D(\varphi^0)$
satisfying $|u_0|_X \le n_2$ and $\varphi^0(u_0)\le n_3$.
\end{description}
} 
\vspace{0.25cm}

\begin{center}
{\bf References}
\end{center}
\begin{enumerate}
\item H.~Attouch, {\it Variational convergence for functions and operators,} Pitman Advanced Publishing Program, {\bf Vol.\,1}, 1984.
\item J.~P.~Aubin, Un th\'eor\`eme de compacit\'e, C.~R.~Acad.~Sci.~Paris, 
{\bf 256} (1963), 5042--5044.
\item V.~Barbu, {\it Nonlinear Differential
 Equations of Monotone Types in Banach Spaces,}
Springer Monographs in Mathematics, Springer, 2010. 
\item H.~Br\'ezis, {\it Op\'eratuers Maximaux Monotones et 
Semi-groupes de 
Contraction dans les Espaces de Hibert}, Math.~Studies {\bf 5}, North-Holland,
Amsterdam, 1973.
\item A.~Damlamian, Some results on the multi-phase Stefan problem, Comm.~Partial Diff.~Eq.\,{\bf 2} (1977), 1017--1044.
\item A.~Damlamian and N.~Kenmochi, Evolution equations generated by
subdifferentials in the dual space of $H^1(\Omega)$,
Discrete Contin. Dynam. Systems {\bf 5} (1999), 269--278.
\item H.~J.~Eberl, M.~Efendiev, D.~Wrzosek and A.~Zhigun, Analysis of a 
degenerate biofilm model with a nutrient taxis term, Discrete Contin. Dynam. Systems {\bf 34}, No.1 (2014), 99--119. 
\item  H.J.~Eberl, D.F.~Parker and M.C.M.~van Loosdrecht, 
A new deterministic
spatio-temporal continuum model for biofilm development, J.~Theoretical 
Medicine\,{\bf 3} (2001), 161--175.
\item T.~Fukao and N.~Kenmochi, Variational inequality for the Navier-Stokes
equations with time-dependent constraint, pp.~87--102 in {\it Computational
Science 2011}, Gakuto Intern.~Math.~Sci.~Appl.
{\bf Vol.\,34}, Gakk\=otosho, Tokyo, 2011.
\item T.~Fukao and N.~Kenmochi, 
Parabolic variational inequalities with weakly time-dependent constraints,  Adv.~Math.~Sci.~Appl., {\bf 23} (2013), 365--395. 
\item T.~Fukao and N.~Kenmochi, Quasi-variational inequalities 
approach to 
heat convection problems with temperature dependent velocity constraint, Discrete Contin. Dynam. Systems 
{\bf 35} (2015), 2523--2538.
\item M.~Gokieli, N.~Kenmochi and M.~Niezg\'odka, Variational inequalities of
Navier-Stokes type with time-dependent constraints, J.~Math.~Anal.~Appl.\,{\bf 449} (2017), 1229--1247.
\item N.~Kenmochi, Solvability of nonlinear evolution equations with 
time-dependent constraints and applications, Bull.~Fac.~Edu., Chiba Univ., {\bf 30} (1981), 1--87.
\item N.~Kenmochi, Monotonicity and compactness methods for nonlinear
variational inequalities, pp.~203--298 in {\it Handbook of Differential 
Equations: Stationary
Partial Differential Equations} {\bf Vol.\,4}, Elsevier, Amsterdam, 2007.
\item N.~Kenmochi and M.~Niezg\'odka, Weak solvability for parabolic 
variational inclusions and application to quasi--variational problems,
Adv.~Math.~Sci.~Appl.\,{\bf 25} (2016), 62--97.
\item U.~Mosco, Convergence of convex sets and of solutions of variational
inequalities, Advances Math.~{\bf 3} (1969), 510--585.
\item  M.~Peszynska, A.~Trykozko, G.~Iltis and S.~Schlueter, 
Biofilm growth in
porous media: Experiments, computational modeling at the porescale, 
and upscaling, Advances in Water Resources (2015), 1--14.
\item K.~Shirakawa, A.~Ito, N.~Yamazaki and N.~Kenmochi, Asymptotic stability
for evolution equations governed by subdifferentials, pp.~287-310 in {\it
Recent Development in Domain Decomposition Methods and Flow Problems}, Gakuto
International Series, Math.~Sci.~Appl., {\bf Vol.\,11}, Gakk\=otosho, Tokyo, 1998.
\end{enumerate}

\end{document}